\documentclass[11pt,a4paper]{article}
\usepackage{amsthm,mathrsfs,amsmath,amssymb,enumerate}
\usepackage{dsfont}
\usepackage{comment}
\usepackage{indentfirst}
\numberwithin{equation}{section}
\makeatother
\usepackage[colorlinks,
linkcolor=blue,
anchorcolor=blue,
citecolor=blue]{hyperref}
\usepackage{hyperref}
\hypersetup{hypertex=true,
	colorlinks=true,
	linkcolor=blue,
	anchorcolor=blue,
	citecolor=blue}
\usepackage{hyperref,amsmath,amssymb,amscd}
\usepackage{amssymb,amsmath,enumerate}
\usepackage{indentfirst}
\usepackage{cite}

\newtheorem{thm}{Theorem}[section]

\newtheorem{lem}[thm]{Lemma}
\newtheorem{prop}[thm]{Proposition}
\theoremstyle{definition}

\newtheorem{rem}[thm]{\bf Remark}

\newenvironment{prf}{\noindent {\bf Proof.}\rm}{\qed}

\allowdisplaybreaks

\begin{document}

\baselineskip 14.5pt \setcounter{page}{1}
\title{\bf Existence and asymptotic properties of standing waves for dipolar Bose-Einstein condensate with rotation\thanks{
Supported by National Natural Science Foundation of China (No. 12371120).}}

\author{{Meng-Hui Wu, Shubin Yu, Chun-Lei Tang\footnote{
Corresponding author. E-mail address: wumenghui7@163.com (M.-H. Wu), yshubin168@163.com (S. Yu), tangcl@swu.edu.cn (C.-L. Tang).} }\\
{\small \emph{ School  of  Mathematics  and  Statistics, Southwest University,  Chongqing {\rm400715},}}\\
{\small \emph{People's Republic of China}}\\}
\date{}
\maketitle

\baselineskip 14.5pt
{{{\bf Abstract}:
In this article, we study the existence and asymptotic properties of prescribed mass standing waves
for the rotating dipolar Gross-Pitaevskii equation with a harmonic potential
in the unstable regime. This equation arises as an effective model describing Bose-Einstein condensate of trapped dipolar quantum gases rotating at the speed $\Omega$.
To be precise, we mainly focus on the two cases:
the rotational speed $0<\Omega<\Omega^{*}$ and $\Omega=\Omega^{*}$, where $\Omega^{*}$ is called a  critical rotational speed.
For the first case, we obtain two different standing waves,
one of which is a local minimizer and can be determined as the ground state, and the other is mountain pass type.
For the critical case,
we rewrite
 the original problem as a dipole Gross-Pitaevskii equation
with a constant magnetic
field and partial harmonic confinement. Under this setting,
 a local minimizer can also be obtained, which seems to be the optimal result.
Particularly, in both cases, we establish the mass collapse behavior of the local minimizers.
Our results
extend the work of Dinh (Lett. Math. Phys., 2022) and
Luo et al. (J. Differ. Equ., 2021) to
the non-axially symmetric harmonic potential, and answer the open
question proposed by Dinh.
}}
\\



\section{Introduction}

In recent years, the dipolar Bose-Einstein condensate (BEC) in a rotating frame have been investigated in
several studies by Physicists \cite{Antoine2018,Tang2017,Zhang2016}.
At temperature much smaller than the critical temperature, the three-dimensional (3D) rotating
dipolar BEC can be well described by the Gross-Pitaevskii equation (GPE) with dipolar interaction \cite{Dell2007,Yi2006} as follows:
\begin{align}\label{EQ1}
i\hbar\frac{\partial \psi(x,t)}{\partial t}=-\frac{\hbar^2}{2m}\nabla^2\psi+V(x)\psi+U_0|\psi|^2\psi
+(V_{dip}\ast|\psi|^2)\psi-\Omega L_{z}\psi,\,\,\,\,x\in\mathbb{R}^3,\,\,\,t>0,
\end{align}
where $x=(x_1,x_2,x_3)^{T}\in\mathbb{R}^3$ is the Cartesian coordinates, $t$ is the time variable,
$\ast$ denotes the convolution operator with respect to the spatial variable, $\hbar$ is the Planck constant,
$m$ is the mass of a dipolar particle, $\Omega\geq 0$ is the angular velocity of the laser beam and
$L_{z}$ is the third component of angular momentum operator $\mathbb L:=-i\hbar x\wedge\nabla$ i.e.,
$L_{z}:=i\hbar(x_2\partial_{x_1}-x_1\partial_{x_2}).$
Moreover, $V(x)$ is the harmonic trapping potential described as
$V(x):=\frac{m}{2}(\gamma_1^2x_1^2+\gamma_2^2x_2^2+\gamma_3^2x_3^2),$
where $\gamma_1$, $\gamma_2$ and $\gamma_3$ being the trap frequencies in each spatial direction. The parameter
$U_0=\frac{4\pi\hbar^2a_{s}}{m}$ represents the short-range interaction
between dipoles in the condensate with $a_{s}$ the $s$-wave scattering length (positive for repulsive interaction and negative for attractive interaction).
The long-range dipolar interaction potential between two dipoles is defined by
\begin{align*}
V_{dip}(x)=\frac{\mu_0\mu_{dip}^2}{4\pi}\frac{1-3\cos^2(\theta)}{|x|^3},\,\,\,x\in\mathbb{R}^3,
\end{align*}
where $\mu_0$ is the vacuum magnetic permeability, $\mu_{dip}$ is the
permanent magnetic dipole moment and $\theta$ is the angle between the vector $x$ and the unitary dipole axis $n$.
In view of the mass conservation, the wave function $\psi(x,t)$ is always normalized such that $\int_{\mathbb{R}^3}|\psi(x,t)|^2dx=N$, where $N$ is the total number of
dipolar particles in the rotating dipolar BEC.

To our knowledge, the rigourous mathematical study of \eqref{EQ1} is only given in \cite[Chapter 4]{CaiYY2011}.
Under suitable
restriction on the rotational speed, they studied the existence and nonexistence of ground states of the 3D rotational
Gross-Pitaevskii-Poisson system (GPPS) and the quasi-2D rotational dipolar
GPE (which are the some dimensionless forms of \eqref{EQ1}),
and computed the ground states for the quasi-2D rotational dipolar
GPE by using a numerical method.

In order to investigate the more mathematical features of the rotating dipolar BEC,
we give the dimensionless form of \eqref{EQ1} by introducing the dimensionless variables as in \cite{Jeanjean2017}, that is,
\begin{align}\label{EQX}
i\partial_{t}\psi=-\frac{1}{2}\Delta \psi+V(x)\psi+\lambda_1|\psi|^2\psi+\lambda_2(K\ast|\psi|^2)\psi-\Omega L_{z}\psi,\,\,\,\,\,\,
(t,x)\in \mathbb{R}^{+}\times\mathbb{R}^3,
\end{align}
where $\Omega\geq0$ is the rotational speed,
the parameters $\lambda_{i}\in\mathbb{R}(i=1,2)$ describe the strength of the two nonlinearities (i.e. the local and the non-local one),
the quantum mechanical angular momentum operator
$$
L_{z}=i(x_2\partial_{x_1}-x_1\partial_{x_2})
$$
and
\begin{align}\label{VX}
V(x)=\frac{1}{2}\sum\limits_{j=1}\limits^{3}\gamma_{j}^2x_{j}^2
\ \mbox{with}\ \gamma_{j}>0
\end{align}
is the dimensionless harmonic trapping potential.
In addition,
the dimensionless long-range dipolar interaction potential is given by
$$
K(x)=\frac{1-3cos^2(\theta)}{|x|^3},\ x\in\mathbb{R}^3,
$$
where $\theta=\theta(x)$
being
the angle between the dipole axis $n$ and the vector $x$. Here to simplify notation,
we fix, without restriction of generality, the dipole axis
$n=(0,0,1).$
For the general dipole axis, we refer to \cite{Bao2010,Bao2012}.
The corresponding normalization is now
$$\int_{\mathbb{R}^3}|\psi(x,t)|^2dx
=\int_{\mathbb{R}^3}|\psi(x,0)|^2dx=c>0.$$

Based on the dimensionless form \eqref{EQX}, in this paper,
we are mainly interesting in the stationary states/standing waves of problem \eqref{EQ1}, that is,
\begin{equation}\label{eqn:psi-definition}
\psi(x,t)=e^{i\mu t}u(x),\,\,\, x\in \mathbb{R}^3,
\end{equation}
where $\mu\in\mathbb{R}$ and $u(x)$ is a time-independent function.
As we will see, we present
new contributions respecting the existence of standing waves for \eqref{EQ1} with $0<\Omega<\Omega^*$ or $\Omega=\Omega^*$, as well as their multiplicity, ground states and asymptotic properties.
Here $\Omega^*>0$ is the critical
rotational speed defined in \eqref{critical-speed} below.
Plugging \eqref{eqn:psi-definition} into \eqref{EQX},  it is natural to consider the following stationary
equation
\begin{align}\label{EQ}
-\frac{1}{2}\Delta u+V(x)u+\lambda_{1}|u|^2u+\lambda_2(K*|u|^2)u-\Omega L_{z}u+\mu u=0\,\,\,\,{\rm{in}}\,\,\mathbb{R}^3,
\end{align}
having the prescribed mass
\begin{equation}\label{eqn:mass}
\int_{\mathbb{R}^3}|u|^2dx=c>0,
\end{equation}
where $\mu\in\mathbb R$ is the Lagrange multiplier.
In particular, we constrain the
coordinate plane $(\lambda_1,\lambda_2)$ to
the so-called unstable regime
\begin{align}\label{UR}
D_{ur}:=\left\{(\lambda_1, \lambda_2)\in\mathbb{R}^2:\lambda_1-\frac{4\pi}{3}\lambda_2<0, \lambda_2>0\,\, \mbox{or}\,\,
  \lambda_1+\frac{8\pi}{3}\lambda_2<0, \lambda_2<0\right\},
\end{align}
which is introduced in \cite[Remark 4.5]{Carles2008} as a  complement
of the
stable regime
\begin{align}\label{SR}
D_{sr}:=\left\{(\lambda_1, \lambda_2)\in\mathbb{R}^2: \lambda_1-\frac{4\pi}{3}\lambda_2\geq0,\lambda_2\geq0\,\, \mbox{or}\,\,
  \lambda_1+\frac{8\pi}{3}\lambda_2\geq0, \lambda_2\leq0\right\}.
\end{align}
Here we remark that
Carles, Markowich and Sparber in \cite{Carles2008} considered
the existence and uniqueness of ground states as well as
the blowup theory for \eqref{EQX} with $\Omega=0$.
Moreover, some crucial results for the dipolar term were established,
such as the Fourier transform of $K(x)$ in $\mathbb{R}^3$ under the condition of $n=(0,0,1)$.
For more related results of \eqref{EQX} and its the stationary form \eqref{EQ} with $\Omega=0$,
we refer to \cite{Dinh2021,Feng2021,Liu2015,Ma2011,Ma2013} for the unstable regime, and \cite{Carles2015} for the stable regime.

In the following, we always assume the rotational speed $\Omega>0$.
In order to expound our main results, we define the working space $\Sigma$ for \eqref{EQ} as follows
\begin{align*}
\Sigma:=\left\{u\in H^1(\mathbb{R}^3, \mathbb{C}):\int_{\mathbb{R}^3}|x|^2|u|^2dx<+\infty\right\},
\end{align*}
which is a Hilbert space with the inner product and norm
\begin{align*}
&(u,v)_{\Sigma}:={\rm{Re}}\int_{\mathbb{R}^3}\left(\nabla u\cdot\nabla\overline{v}+|x|^2u\overline{v}+u\overline{v}\right)dx,\\
&\|u\|_{\Sigma}:=\left(\int_{\mathbb{R}^3}|\nabla u|^2+|x|^2|u|^2+|u|^2dx\right)^{\frac{1}{2}}.
\end{align*}
Here ``Re" stays for the real part and $\overline{v}$ denotes the conjugate of $v$. Moreover, we define
$$
\|u\|_{\dot{\Sigma}}:=\left(\int_{\mathbb{R}^3}|\nabla u|^2+|x|^2|u|^2dx\right)^{\frac{1}{2}},
$$
which is a equivalent norm of $\Sigma$.
To obtain the existence of solutions (also known as normalized solutions) for \eqref{EQ}, we need to seek
critical points of the energy functional
\begin{align}\label{FH}
E_{\Omega}(u)\!:=\frac{1}{2}\int_{\mathbb{R}^3}\!|\nabla u|^2\!+V(x)|u|^2dx
\!+\frac{1}{2}\int_{\mathbb{R}^3}\lambda_1|u|^4
\!+\lambda_2(K\ast|u|^2)|u|^2dx
\!-\Omega\int_{\mathbb{R}^3}\overline{u}L_{z}udx
\end{align}
under the mass constraint
\begin{align*}
S(c):=\{u\in \Sigma:|u|_2^2=c\}.
\end{align*}
By the standard argument, we know that $E_{\Omega}(u)$ is well-defined on $\Sigma$ (see \eqref{E1}).
%
%
Furthermore, it is well-known that for the harmonic trapping potential in \eqref{VX},
there is a
critical rotational speed
\begin{align}\label{critical-speed}
\Omega^{*}:=\min\{\gamma_1, \gamma_2\}.
\end{align}
If $\Omega<\Omega^*$, we can get the equivalence of a new norm involving the rotational term,
that is,
\begin{align*}
\int_{\mathbb{R}^3}|\nabla u|^2+2V(x)|u|^2-2\Omega\overline{u} L_{z}udx\cong\int_{\mathbb{R}^3}|\nabla u|^2+|x|^2|u|^2dx,
\end{align*}
see Lemma \ref{lem5}.
In this sense, one can overcome the difficulty (i.e. how to recover the compactness) brought by the rotational term. Indeed, this equivalence has been observed in
\cite{Arbunich2019,Dinh2022,Luo2021,Lewin2018}, where the authors considered the classical
rotating BEC in the mean-field regime described by the nonlinear GPE with pure power nonlinearities.
Compared to the classical
rotating BEC, our problem involves the nonlinear
potential energy
\begin{align}\label{R}
N(u):=\int_{\mathbb{R}^3}\lambda_1|u|^4+\lambda_2(K*|u|^2)|u|^2dx,
\end{align}
which is sign indefinite if $(\lambda_1, \lambda_2)\in D_{ur}$ but there is a $u\in S(c)\cap C_0^{\infty}(\mathbb{R}^3)$ such that $N(u)<0$ (see \cite{Jeanjean2017}).
From this perspective, the unstable regime is very interesting.
In particular, for the stable regime \eqref{SR},
  we know that $N(u)>0$. In this case,
 we refer to \cite[Theorem 1]{Triay2018} about the existence of ground states of \eqref{EQ} if $\Omega<\Omega^*$.

Now we give the first main result for the existence of normalized solutions to problem \eqref{EQ} if $\Omega<\Omega^*$ and $(\lambda_1, \lambda_2)\in D_{ur}$.
This result implies that  problem \eqref{EQ} has
 at least two solutions, which correspond to a local minimizer and a mountain pass critical point of $E_{\Omega}|_{S(c)}$, respectively.

\begin{thm}\label{thm1}
Let $V$ be as in \eqref{VX}, $0<\Omega<\Omega^{*}=\min\{\gamma_1, \gamma_2\}$ and $(\lambda_1, \lambda_2)\in D_{ur}$. Then
\begin{enumerate}[\indent$(i)$]
  \item for any fixed $r>0$,
there exists $c_0:=c_0(r,\Omega)>0$ {\rm(}see \eqref{CS}{\rm)} such that for any $0<c<c_0$, there is an interior local minimizer $u_{c}$ of $E_{\Omega}(u)$ on the set $S(c)\cap B(r)$, where $B(r):=\{u\in\Sigma: \|u\|_{\dot{\Sigma}}^2\leq r\}.$
Moreover, $u_c$ is a weak solution of problem \eqref{EQ} with the Lagrange multiplier $\mu_c<0$ and
\begin{align}\label{mu-c-bound}
\mu_c\in\left[-\zeta^0-\delta(r)c,3\left(-\frac{C_{*}}{2}+\Lambda \mathcal C_4^4r^{\frac{1}{2}}c^{\frac{1}{2}}\right)\right],
\end{align}
where $\zeta^0:=\frac{1}{2}\sum\limits_{j=1}\limits^{3}\gamma_{j}$, $\delta(r):=\sqrt{2}\Lambda \mathcal{C}_4^4(\zeta^0)^{\frac{3}{2}}+\frac{\Lambda}{2}\mathcal C_4^4r^{\frac{3}{2}}$, $C_{*}:=\min\{\frac{(\Omega^{*})^2-\Omega^2}{(\Omega^{*})^2+\Omega^2}, \frac{(\Omega^{*})^2-\Omega^2}{2}, \gamma_3^2\}$, and $\mathcal C_4$, $\Lambda$ are positive constants.
  \item If $0<\Omega<\frac{\sqrt{5}}{3}\Omega^{*}$, then for any $0<c<c_0$,
$E_{\Omega}|_{S(c)}$ has a
second critical point of mountain pass type $v_{c}$ at the level $\gamma(c)> E_{\Omega}(u_{c})$. Moreover,
$v_{c}$ is also a weak solution of problem \eqref{EQ} with the Lagrange multiplier $\widehat{\mu}_{c}\in\mathbb R$.
\end{enumerate}
\end{thm}

\begin{rem}\label{mountain-pass}
Note that we restrict $0<\Omega<\frac{\sqrt{5}}{3}\Omega^{*}$ when finding the mountain pass solution.
We mention that the condition is  pivotal for proving the boundedness of a specific Palai-Smale sequence $\{v_{n}\}$ constructed by a variant of the min-max theorem in \cite{Jeanjean1997}.
Indeed,  it follows from \cite{Luo2021} that the Pohozaev identity \eqref{D5} does not involve the rotation term .
Thus, due to the fact that the rotation term $\int_{\mathbb{R}^3}\Omega\overline{v}_{n}L_{z}v_{n}dx$ within $E_{\Omega}(v_{n})$ is sign indefinite, we  introduce this stronger condition to ensure the boundedness of $\{v_{n}\}$ in $\Sigma$,
see \eqref{bounded-PS-sequence}. Once the boundedness is established, the compactness is natural since $\Sigma\hookrightarrow L^{p}(\mathbb{R}^3)$ is compact for $p\in[2,6)$.
\end{rem}

Next we provide some asymptotic properties of
the local minimizer $u_{c}$ obtained by Theorem \ref{thm1} as $c\rightarrow0^+$.
To this end, we define
\begin{align*}
\mathcal{M}_{c}^{r}:=\{u\in S(c)\cap B(r): E_{\Omega}(u)=m(c,r)\},
\end{align*}
where $m(c,r)=\inf_{u\in S(c)\cap B(r)} E_{\Omega}(u)$ is the minimal energy.
Based on Theorem \ref{thm1}, we know that
$\mathcal{M}_{c}^{r}\not=\emptyset.$
In particular, we will prove that $u_{c}$ is a ground state if  $0<\Omega<\frac{\sqrt{5}}{3}\Omega^{*}$ and $c>0$ sufficiently small. The definition of ground states can be seen in \cite[Definition 1.1]{Jeanjean2017}.
\begin{thm}\label{thm1''}
Under the setting of Theorem {\rm\ref{thm1}}, there holds
\begin{align}\label{11}
\sup\limits_{u\in \mathcal{M}_{c}^{r}}\|u-\varrho_0\psi_0\|^2_{\dot{\Sigma}}=O(c+c^2),
\end{align}
where $\psi_0$ is the unique normalized positive eigenfunction
of harmonic oscillator $-\frac{1}{2}\Delta+\frac{|x|^2}{2}$ and $\varrho_0:=\int_{\mathbb{R}^3}u\psi_0dx$.
In particular, if $\zeta^0=\frac{3C_{*}}{2}$, then
\begin{align}\label{11}
\sup\limits_{u\in \mathcal{M}_{c}^{r}}\|u-\varrho_0\psi_0\|^2_{\dot{\Sigma}}=O(c^2).
\end{align}
\end{thm}
\begin{rem}
The above result implies that the local minimizers in $\mathcal{M}_{c}^{r}$  behave like the first eigenfunction of the harmonic oscillator for small $c>0$. Moreover, \eqref{11} gives the more refined estimate, which
strictly depends on $\zeta^0=\frac{3C_{*}}{2}$.
Recall that $\zeta^0=\frac{1}{2}\sum\limits_{j=1}\limits^{3}\gamma_{j}$
and $C_{*}=\min\{\frac{(\Omega^{*})^2-\Omega^2}{(\Omega^{*})^2+\Omega^2}, \frac{(\Omega^{*})^2-\Omega^2}{2}, \gamma_3^2\}$.
If $\gamma_1=\gamma_2$, we can see that
$\zeta^0 <\frac{3C_{*}}{2}$.
For $\gamma_1\not=\gamma_2$,
it is easy to check that the assumption $\zeta^0=\frac{3C_{*}}{2}$ is well-defined when $\Omega>0$ approaches $0$ and $C_{*}= \frac{(\Omega^{*})^2-\Omega^2}{2}$. Thus,
in this sense, this estimate extends the result in \cite[Theorem 1.1]{Luo2021}.
\end{rem}

\begin{thm}\label{thm1'}
Let $(u_{c}, \mu_{c})\in\mathcal{M}_{c}^{r}\times\mathbb{R}$ be given by Theorem {\rm\ref{thm1}}.
Then $u_{c}\to0$ in $\Sigma$ as $c\to0^{+}$,
$$\lim\limits_{c\to0^{+}}\frac{m(c,r)}{c}
=\lim\limits_{c\to0^{+}}-\mu_{c}=-\mu_0$$ for some $\mu_0\in[-\zeta^0, -\frac{3C_{*}}{2}]$, and
\begin{align*}
\lim\limits_{c\to0^{+}}\frac{\int_{\mathbb{R}^3}|\nabla u_{c}|^2dx-\Omega\int_{\mathbb{R}^3}\overline{u}_{c}L_{z}u_{c}dx}{c}
=\lim\limits_{c\to0^{+}}\frac{\int_{\mathbb{R}^3}2V(x)|u_{c}|_2^2-\Omega\int_{\mathbb{R}^3}\overline{u}_{c}L_{z}u_{c}dx}{c}=-\mu_0.
\end{align*}
Moreover, if $0<\Omega<\frac{\sqrt{5}}{3}\Omega^{*}$, $u_{c}$ is a normalized ground state to \eqref{EQ}-\eqref{eqn:mass} when $c>0$ is sufficiently small.
\end{thm}
\begin{rem}
Theorem \ref{thm1'}
describes a mass collapse behavior of the local minimizers.
Moreover, it can be seen that the local minimizer $u_{c}$ is a ground state when $c$ approaches $0$ and $\Omega>0$ is small.
Indeed,
in order to determine this fact,
a direct way is to prove that all solutions with energy strictly less than $m(c,r)$ are in $S(c)\cap B(r)$ provided that $c\rightarrow0^+$.
Based on this, the Pohozaev identity will also play a crucial role, but similar phenomenon as in Remark \ref{mountain-pass}  will arise. Therefore, we also restrict that $0<\Omega<\frac{\sqrt{5}}{3}\Omega^{*}$.
\end{rem}

As mentioned earlier, there is a critical
rotational speed $\Omega^*$ and we establish the existence and asymptotic properties of
standing waves if $\Omega<\Omega^*$ (i.e. the low rotational speed). Hence, a natural question is whether we can obtain similar results if $\Omega=\Omega^*$ (i.e. the critical rotational speed).
Indeed, the situation is even more subtle. Since the norm equivalence is no longer available, the working space $ \Sigma$ is not helpful in finding normalized solutions to problem
\eqref{EQ} with $\Omega=\Omega^*$.
It leads us to seek new method to deal with the case $0<\Omega=\Omega^{*}$.
Motivated by \cite{Dinh2022},
we rewrite \eqref{EQ} as the following  dipolar GPE with magnetic  field
\begin{align}\label{EQ-magnetic}
-\frac{1}{2}(\nabla-iA)^2u+V_{\Omega^*}(x)u+\lambda_{1}|u|^2u
+\lambda_2(K*|u|^2)u+\mu u=0\,\,\,\,{\rm{in}}\,\,\mathbb{R}^3,
\end{align}
where $A=A(x):=\Omega^*(-x_2, x_1,0)$, $V_{\Omega^*}(x):=V(x)-\frac{1}{2}|A|^2$ and the magnetic Schr\"{o}dinger operator (see \cite{Esteban1989}) is defined by
$$
-(\nabla-iA)^2u:=-\Delta u+2i A\cdot \nabla u-iu{\rm div} A+|A|^2u
=-\Delta u+2i A\cdot \nabla u+|A|^2u.
$$
Recall that $\Omega^{*}=\min\{\gamma_1,\gamma_2\}$.
Without loss of generality, we assume that
$$
\min\{\gamma_1,\gamma_2\}=\gamma_2.
$$
Then $\Omega^{*}=\gamma_2$ and $$
V_{\Omega^{*}}(x)=\frac{1}{2}(\gamma_1^2-\gamma_2^2)
x_1^2+\frac{1}{2}\gamma_3^2x_3^2.$$
From this perspective, \eqref{EQ-magnetic} can be viewed as the dipolar GPE with a constant magnetic field and a partial harmonic confinement in $x_1$, $x_3$ directions.
Moreover, the corresponding energy functional $E_{\Omega^{*}}(u)$ can be rewritten as
\begin{align}\label{4.1}
E_{\Omega^{*}}(u)=
\frac{1}{2}\int_{\mathbb{R}^3}|(\nabla-iA)u|^2dx
+\int_{\mathbb{R}^3}V_{\Omega^*}(x)|u|^2dx
+\frac{1}{2}N(u),
\end{align}
where $N(u)$ is shown in \eqref{R}.

In particular, the working space  is naturally replaced with $\Sigma_{\Omega^{*}}$,
 where
\begin{align}\label{4.2}
\Sigma_{\Omega^{*}}:=\{u\in H_{A}^1(\mathbb{R}^3, \mathbb{C}):\int_{\mathbb{R}^3}V_{\Omega^{*}}(x)|u|^2dx<\infty\}
\end{align}
is equipped with the norm
\begin{align}\label{4.3}
\|u\|_{\Sigma_{\Omega^{*}}}^2:=\int_{\mathbb{R}^3}|(\nabla-iA)u|^2
+2V_{\Omega^*}(x)|u|^2+|u|^2dx
\end{align}
and $H_{A}^1(\mathbb{R}^3, \mathbb{C})$ is the usual magnetic Sobolev space defined by
\begin{align*}
H_{A}^1(\mathbb{R}^3, \mathbb{C}):=\{u\in L^2(\mathbb{R}^3):|(\nabla-iA)u|\in L^2(\mathbb{R}^3)\}.
\end{align*}
Obviously, the energy functional $E_{\Omega^*}$ is well-defined on $\Sigma_{\Omega^{*}}$.
Some basic properties of the magnetic Sobolev space $H_{A}^1(\mathbb{R}^3, \mathbb{C})$ are collected in
\cite[Lemma 2.3]{Dinh2022}. Here we only mention the  diamagnetic inequality
\begin{align}\label{DI}
|\nabla|u(x)||\leq|(\nabla-iA)u(x)|, \,\,\,\,\,a.e.\,\,\,x\in\mathbb{R}^3,
\end{align}
see \cite{Lieb2001}.
On the other hand, unlike $\Omega<\Omega^*$, we further constrain
$(\lambda_1, \lambda_2)$ to belong to a special unstable region $\bar D_{ur}$, where
\begin{align}\label{UR'}
\bar D_{ur}:=\left\{(\lambda_1, \lambda_2)\in\mathbb{R}^2:
\lambda_1+\frac{8\pi}{3}\lambda_2<0, \lambda_2>0,\ \mbox{or}\
  \lambda_1-\frac{4\pi}{3}\lambda_2<0, \lambda_2<0\right\}.
\end{align}
If $(\lambda_1, \lambda_2)\in \bar D_{ur}$, it is esay to check that $N(u)<0$ for all $u\in \Sigma_{\Omega^{*}}$. As we will see, this property plays some important roles, such as excluding the vanishing case for the minimizing sequences of the local minimization problem
\begin{align*}
m_{\Omega^{*}}(c,\rho)=\inf\{E_{\Omega^{*}}(u): u\in S(c)\cap B_{\Omega^{*}}(\rho)\},
\end{align*}
where
$B_{\Omega^{*}}(\rho)$ defined in \eqref{B-rho-definition} and $S(c)=\{u\in \Sigma_{\Omega^*}:|u|_2^2=c\}$.
Furthermore, we point out that for $\Omega=\Omega^*$, inspired by the idea of \cite{Bellazzini2017}, an interior local minimizer can be obtained for the minimization problem \eqref{4.4}. However,
since $\Omega=\Omega^*$ is considered,
we cannot restrict $\Omega$ to be small and then the Pohozaev identity does not seem to provide us any assistance.
Thus, the existence of ground states or mountain pass solutions cannot be determined.
In this sense, our result seems to be optimal.
Especially, we mention that
the authors in \cite{Bellazzini2017} studied the nonlinear Schr\"{o}dinger equation  with a partial confinement, the ground state can be acquired, which is different from our problem.

According to the above framework, we now present the results for the existence and some asymptotic properties of normalized solutions for \eqref{EQ} with $0<\Omega=\Omega^{*}$.

\begin{thm}\label{thm3}
Let $V$ be as in \eqref{VX} satisfying $\gamma_2\leq\gamma_1
<\sqrt{\frac{2M(\gamma_2,\gamma_3)\gamma_2}{c}+\gamma_2^2}$,
$(\lambda_1, \lambda_2)\in \bar D_{ur}$ and $0<\Omega=\Omega^{*}$ with $\Omega^{*}=\gamma_2$,
 where $M(\gamma_2, \gamma_3)$ is a constant depending on $\gamma_2$, $\gamma_3$. Then
\begin{enumerate}[\indent$(i)$]
 \item for any fixed $\rho>0$, there exists $\bar c(\rho)>0$ {\rm(}see \eqref{bar-c-rho}{\rm)}
such that for all $0<c<\bar c(\rho)$, $E_{\Omega^{*}}(u)$ has an interior local minimizer $\widetilde{u}_{c}$ such that
\begin{align*}
E_{\Omega^{*}}(\widetilde{u}_{c})=m_{\Omega^{*}}(c,\rho)<\gamma c.
\end{align*}
Moreover, $\widetilde{u}_{c}$ is a weak solution of problem \eqref{EQ}
with the Lagrange multiplier $\widetilde{\mu}_{c}<0$ and
\begin{align}\label{T3}
\widetilde{\mu}_{c}\in
\Big(\!\!-\!\gamma,\Big(-\frac{1}{2}+\Lambda \mathcal C_4^4\rho^{\frac{1}{2}}c^{\frac{1}{2}}\Big)\gamma
\Big],
\end{align}
where $\Lambda>0$ is a constant and $\gamma:=\gamma_2+\frac{\gamma_3}{2}$.

\item
$\widetilde{u}_{c}\to0$ in $\Sigma_{\Omega^{*}}$ as $c\to0^{+}$, and
$$\lim\limits_{c\to0^{+}}\frac{m_{\Omega^{*}}
(c,\rho)}{c}=\lim\limits_{c\to0^{+}}-\widetilde{\mu}_{c}
=-\widetilde{\mu}_0$$ for some $\widetilde{\mu}_0\in[-\gamma,-\frac{1}{2}\gamma]$.
\end{enumerate}
\end{thm}
\begin{rem}
Compared with \cite[Theorem 1.4]{Dinh2022},
 Theorem \ref{thm3} shows the existence of normalized solutions for \eqref{EQ}
with $\Omega=\Omega^{*}$ when the harmonic potential is not axially symmetric, i.e., $\gamma_1\neq\gamma_2$.
This solves the open problem
proposed in \cite[Remark 1.7]{Dinh2022}.
\end{rem}

\begin{rem}
To sum up, we obtain the existence of prescribed mass/normalized  standing waves for \eqref{EQX} with $0<\Omega\leq \Omega^*$.  For the stability,
 we believe that corresponding result can be established based on the arguments in \cite{Jeanjean2017,Dinh2022}, and this will be treated as an open question.
Finally, we remark that for $\Omega>\Omega^{*}$,
some non-existence results can be seen in \cite{CaiYY2011,Bao2005} and the numerical method was used.
\end{rem}

The paper is organized as follows. In Section \ref{Preliminaries}, some preliminary results are given.
Section  \ref{low} is devoted to investigating the case of $0<\Omega<\Omega^*$, and accomplishing the proof of
 Theorems \ref{thm1}-\ref{thm1'}. In Section \ref{critical}, we consider the case of $\Omega=\Omega^*$  and prove Theorem \ref{thm3}.

\section{Preliminaries}\label{Preliminaries}
In this section, we give some preliminary results. In the process of proof, we need to use the following notation.
For $1\leq p<\infty$ and $u\in L^{p}(\mathbb{R}^3,\mathbb{C})$, the norm is denoted by $|u|_{p}:=\left(\int_{\mathbb{R}^3}|u|^{p}dx\right)^{\frac{1}{p}}$.
The Hilbert spaces $H^{1}(\mathbb{R}^3, \mathbb{C})$ is defined as
$$
H^1(\mathbb{R}^3, \mathbb{C}):=\left\{u\in L^2(\mathbb{R}^3, \mathbb{C}): \nabla u\in L^2(\mathbb{R}^3, \mathbb{C})\right\},
$$
with the inner product $(u,v):={\rm{Re}}\left[\int_{\mathbb{R}^3}\nabla u\cdot\nabla \overline{v}+u\overline{v}dx\right]$
and norm $\|u\|:=\left(|\nabla u|_2^2+|u|_2^2\right)^{\frac{1}{2}}$.
We simply write $H^1(\mathbb{R}^3)$ for $H^{1}(\mathbb{R}^3, \mathbb{C})$, and $L^{p}(\mathbb{R}^3)$ for $L^{p}(\mathbb{R}^3,\mathbb{C})$.
$H^{-1}(\mathbb{R}^3)$ is the dual space of $H^1(\mathbb{R}^3)$. $C$ and $C_{i}$($i=1,2,...$) represent different positive constants.
$o_{n}(1)$ denotes a quantity which goes to zero. ``Re" stays for the real part and $\overline{v}$ denotes the conjugate of $v$.

For any $u\in H^1(\mathbb{R}^3)$ and $q\in(2,6)$,  the Gagliardo-Nirenberg inequality \cite{Weinstein2009} as follows: there is a constant $\mathcal{C}_{q}>0$ such that
\begin{align}\label{GN}
|u|_{q}\leq \mathcal{C}_{q}|\nabla u|_2^{\delta_{q}}|u|_2^{(1-\delta_{q})},
\end{align}
where $\delta_{q}$ is defined by $\delta_{q}:=\frac{3(q-2)}{2q}$.

In particular, the proof of the main results need to use some spectral theory. 
Here we introduce
\begin{align}\label{T}
\zeta^0:=\inf\left\{\frac{1}{2}\int_{\mathbb{R}^3}|\nabla u|^2dx+\int_{\mathbb{R}^3}V(x)|u|^2dx: u\in\Sigma,\,\,\,|u|_2^2=1\right\}
\end{align}
as the simple first eigenvalue of the multi-dimensional harmonic oscillator $-\frac{1}{2}\Delta +V$,
where $$V(x)=\frac{1}{2}\sum\limits_{j=1}\limits^{3}\gamma_{j}^2x_{j}^2.$$ Let $\Phi_0(x)\in S(1):=\{u\in\Sigma:|u|_2^2=1\}$ be an eigenfunction with respect to the first
eigenvalue $\zeta^0$. It is well-known that $\zeta^0=\frac{1}{2}\sum\limits_{j=1}\limits^{3}\gamma_{j}$, $\Phi_0(x)$ is real-valued and
\begin{align}\label{E}
\Phi_0(x):=\pi^{-\frac{3}{4}}\left(\prod\limits_{j=1}\limits^{3}\sqrt{\gamma_{j}}\right)^{\frac{1}{2}}e^{-\frac{1}{2}\sum\limits_{j=1}\limits^{3}\gamma_{j}x_{j}^2}.
\end{align}

Let the Fourier transform on the Schwartz space of $u$ by $\widehat{u}(\xi):=\mathcal{F}u(\xi)=\int_{\mathbb{R}^3}e^{-ix\cdot\xi}u(x)dx$.
Then by \cite[Lemma 2.3]{Carles2008}, we have the following result about the term $N(u)$ given in \eqref{R}.

\begin{lem}\label{lem2}
Let $u\in H^1(\mathbb{R}^3)$, then
$$
N(u)=\frac{1}{(2\pi)^3}\int_{\mathbb{R}^3}[\lambda_1+\lambda_2\widehat{K}(\xi)]|\widehat{u^2}|^2d\xi
$$
and $|N(u)|\leq \Lambda|u|_4^4$,
where $\Lambda:=\max\left\{|\lambda_1-\frac{4\pi}{3}\lambda_2|, |\lambda_1+\frac{8\pi}{3}\lambda_2|\right\}$.
\end{lem}

Due to the definition of the potential $V(x)$ in \eqref{VX}, we have the following compactness result, which is verified in \cite[Lemma 3.1]{Zhang2000}.
\begin{lem}\label{lem4}
The embedding $\Sigma\hookrightarrow L^{p}(\mathbb{R}^3)$ is compact for $p\in[2,6)$.
\end{lem}

To discuss the existence of mountain pass solutions to \eqref{EQ}-\eqref{eqn:mass}, we give the following Pohozaev identity,
which can be obtained by the similar arguments in \cite[Proposition 3.2]{Luo2021} and \cite[Lemma 2.2]{Antonelli2011}. So we omit it.
\begin{lem}\label{lem6}
Let $(\lambda_1, \lambda_2)\in D_{ur}$, $\mu\in\mathbb{R}$,
and $0<\Omega<\Omega^{*}$. If $u\in \Sigma$ is a weak solution for the equation
\begin{align}\label{D4}
-\frac{1}{2}\Delta u+V(x)u+\lambda_{1}|u|^2u+\lambda_2(K*|u|^2)u-\Omega L_{z}u+\mu u=0\,\,\,\,{\rm{in}}\,\,\mathbb{R}^3,
\end{align}
where $V(x)$ be as in \eqref{VX},
then the corresponding Pohozaev identity as follows:
\begin{align}\label{D5}
Q(u):=\frac{1}{2}\int_{\mathbb{R}^3}|\nabla u|^2dx-\int_{\mathbb{R}^3}V(x)|u|^2dx+\frac{3}{4}\int_{\mathbb{R}^3}\lambda_1|u|^4+\lambda_2(K\ast|u|^2)|u|^2dx=0.
\end{align}
\end{lem}

Moreover, as stated in introduction, we establish the following equivalent norm for $\Omega<\Omega^*$.
\begin{lem}\label{lem5}
If $0<\Omega<\Omega^{*}$, then for any $u\in \Sigma$,
\begin{align*}
\int_{\mathbb{R}^3}|\nabla u|^2+2V(x)|u|^2-2\Omega\overline{u} L_{z}udx\cong\int_{\mathbb{R}^3}|\nabla u|^2+|x|^2|u|^2dx,
\end{align*}
where $V(x)$ be as in \eqref{VX}.
\end{lem}
\begin{prf}
It follows from H\"{o}lder's, Cauchy-Schwarz's and Young's inequalities that
\begin{align}\label{E1}
&\ \ \ \left|\int_{\mathbb{R}^3}\Omega\overline{u} L_{z}udx\right|\notag\\
&=\left|\Omega\int_{\mathbb{R}^3}i\overline{u}(x_2\partial_{x_1}u-x_1\partial_{x_2}u)dx\right|\notag\\
&\leq
\Omega\left[\left(\int_{\mathbb{R}^3}x_2^2|u|^2dx\right)^{\frac{1}{2}}\left(\int_{\mathbb{R}^3}|\partial_{x_1}u|^2dx\right)^{\frac{1}{2}}+
\left(\int_{\mathbb{R}^3}x_1^2|u|^2dx\right)^{\frac{1}{2}}\left(\int_{\mathbb{R}^3}|\partial_{x_2}u|^2dx\right)^{\frac{1}{2}}\right]\notag\\
&\leq\left(\int_{\mathbb{R}^3}|\nabla u|^2dx\right)^{\frac{1}{2}}\left(\int_{\mathbb{R}^3}\Omega^2(x_1^2+x_2^2)|u|^2dx\right)^{\frac{1}{2}}\notag\\
&\leq\varepsilon\int_{\mathbb{R}^3}|\nabla u|^2dx+\frac{\Omega^2}{4\varepsilon}\int_{\mathbb{R}^3}
(x_1^2+x_2^2)|u|^2dx
\end{align}
for any $\varepsilon>0$. By \eqref{E1}, we have
\begin{align*}
&\int_{\mathbb{R}^3}|\nabla u|^2+2V(x)|u|^2-2\Omega\overline{u} L_{z}udx\\
&\quad\leq(1+2\varepsilon)\int_{\mathbb{R}^3}|\nabla u|^2dx+\int_{\mathbb{R}^3}2V(x)|u|^2dx+\frac{\Omega^2}{2\varepsilon}\int_{\mathbb{R}^3}(x_1^2+x_2^2)|u|^2dx\\
&\quad \leq\max\{1+2\varepsilon, \gamma_1^2+\frac{\Omega^2}{2\varepsilon}, \gamma_2^2+\frac{\Omega^2}{2\varepsilon}, \gamma_3^2\}\int_{\mathbb{R}^3}|\nabla u|^2+|x|^2|u|^2dx\\
&\quad\triangleq C^{*}\int_{\mathbb{R}^3}|\nabla u|^2+|x|^2|u|^2dx.
\end{align*}
To claim the reverse inequality, we apply \eqref{E1} again to deal with the rotation term. So
\begin{align*}
&\int_{\mathbb{R}^3}|\nabla u|^2+2V(x)|u|^2-2\Omega\overline{u} L_{z}udx\\
&\quad\geq(1-2\varepsilon)\int_{\mathbb{R}^3}|\nabla u|^2dx+\int_{\mathbb{R}^3}2V(x)|u|^2dx-\frac{\Omega^2}{2\varepsilon}\int_{\mathbb{R}^3}(x_1^2+x_2^2)|u|^2dx\\
&\quad=(1-2\varepsilon)\int_{\mathbb{R}^3}|\nabla u|^2dx
+\int_{\mathbb{R}^3}((\gamma_1^2-\frac{\Omega^2}{2\varepsilon})x_1^2+(\gamma_2^2-\frac{\Omega^2}{2\varepsilon})x_2^2+\gamma_3^2x_3^2)|u|^2dx.
\end{align*}
In view of $0<\Omega<\Omega^{*}=\min\{\gamma_1, \gamma_2\}$ and the arbitrariness of $\varepsilon$, we choose $\varepsilon=\frac{\Omega^2}{(\Omega^{*})^2+\Omega^2}$ such that
\begin{align*}
1-2\varepsilon=\frac{(\Omega^{*})^2-\Omega^2}{(\Omega^{*})
^2+\Omega^2}>0\quad\mbox{and}\
\gamma_{j}^2-\frac{\Omega^2}{2\varepsilon}\geq(\Omega^{*})^2-\frac{\Omega^2}{2\varepsilon}=\frac{(\Omega^{*})^2-\Omega^2}{2}>0,\,\,\, j=1,2.
\end{align*}
Thus,
\begin{align*}
\int_{\mathbb{R}^3}|\nabla u|^2+2V(x)|u|^2-2\Omega\overline{u} L_{z}udx\geq C_{*}\int_{\mathbb{R}^3}|\nabla u|^2+|x|^2|u|^2dx,
\end{align*}
where $C_{*}:=\min\{\frac{(\Omega^{*})^2-\Omega^2}{(\Omega^{*})^2+\Omega^2}, \frac{(\Omega^{*})^2-\Omega^2}{2}, \gamma_3^2\}$.
Consequently,
$$
C_{*}\int_{\mathbb{R}^3}|\nabla u|^2+|x|^2|u|^2dx\leq\int_{\mathbb{R}^3}|\nabla u|^2+2V(x)|u|^2-2\Omega\overline{u} L_{z}udx\leq
C^{*}\int_{\mathbb{R}^3}|\nabla u|^2+|x|^2|u|^2dx.
$$
This completes the proof.
\end{prf}

\section{The low rotational speed}\label{low}

In this section, we consider the case of $0<\Omega<\Omega^{*}$ (the low rotational speed), and give the proof of Theorems \ref{thm1}, \ref{thm1''} and \ref{thm1'}.

 Note that $(\lambda_1, \lambda_2)\in D_{ur}$, then there exists $u\in S(c)\cap C_0^{\infty}(\mathbb{R}^3)$ such that $N(u)<0$, see \cite{Jeanjean2017}.
Moreover, we have $u_{t}(x):=t^{\frac{3}{2}}u(tx)\in S(c)$ for any $t>0$
and
\begin{align*}
E_{\Omega}(u_{t})=\frac{t^2}{2}\int_{\mathbb{R}^3}|\nabla u|^2dx+\frac{1}{t^2}\int_{\mathbb{R}^3}V(x)|u|^2dx+\frac{t^3}{2}N(u)
-\Omega\int_{\mathbb{R}^3}\overline{u} L_{z}udx\to-\infty
\end{align*}
as $t\to+\infty$, which imply that the usual global minimization method does not work.
But due to the appearance of the harmonic potential $V(x)$, it seems to find
a critical point of $E_{\Omega}\mid_{ S(c)}$ by minimizing on a constructed sub-manifold of $S(c)$. In addition, a mountain pass type critical point will appear.

\subsection{Existence of local minimizers}

In this subsection, inspired by \cite{Jeanjean2017,Luo2021},
we introduce a local minimization problem, that is, for any
given $r>0$,
\begin{align}\label{3.1}
m(c,r):=\inf\limits_{u\in S(c)\cap B(r)} E_{\Omega}(u),
\end{align}
where
$$
B(r):=\{u\in\Sigma: \int_{\mathbb{R}^3}|\nabla u|^2+|x|^2|u|^2dx\leq r\}.
$$
For any fixed $r>0$, it is clear that $m(c,r)>-\infty$ if $S(c)\cap B(r)\neq \emptyset$.
Then we will prove that $S(c)\cap B(r)\neq \emptyset$ (see Lemma \ref{lem6}) and $m(c,r)$ is achieved by applying the compactness property (Lemma \ref{lem4}).
To conclude that this minimizer is indeed a critical point of $E_{\Omega}|_{S(c)}$,
it suffices to exclude the possibility that this minimizer belong to the boundary of $S(c)\cap B(r)$. The corresponding procedures are stated as follows.

First,  the proof of the set $S(c)\cap B(r)$ is not empty is given, which also implies that $m(c,r)$ is well-defined.
\begin{lem}\label{lem6}
Let $0<\Omega<\Omega^{*}$. Then for any $r>0$, there
exists $r_1>0$ such that for all $0<c<\frac{r}{r_1}$,
$$
S(c)\cap B(r)\neq\emptyset.
$$
\end{lem}
\begin{prf}
By \eqref{E}, we define $\Phi_c(x):=\sqrt{c}\Phi_0(x)\in S(c)$ and
 $r_1:=\int_{\mathbb{R}^3}|\nabla \Phi_0|^2+|x|^2|\Phi_0|^2dx$.
  It is obvious that
$\int_{\mathbb{R}^3}|\nabla \Phi_{c}|^2+|x|^2|\Phi_{c}|^2dx=cr_1<r$
for all $0<c<\frac{r}{r_1}$. Thus, we have $\Phi_{c}(x)\in S(c)\cap B(r)$, and this completes the proof.
\end{prf}

Now, we prove the existence of local minimizers for $m(c,r)$  and obtain Theorem \ref{thm1}-$(i)$.

\noindent\textbf{Proof of Theorem \ref{thm1} (existence of local minimizers)}.
By Lemma \ref{lem6},
we know that for any given $r>0$, $m(c,r)>-\infty$ for all $0<c<\frac{r}{r_1}$. Let $\{u_{n}\}\subseteq S(c)\cap B(r)$ be a minimizing
sequence for $m(c,r)$ satisfying $\lim\limits_{n\to\infty} E_{\Omega}(u_{n})=m(c,r)$,
then $\{u_{n}\}$ is bounded in $\Sigma$. By Lemma \ref{lem4}, there exists $u_{c}\in\Sigma$
such that $u_{n}\rightharpoonup u_{c}$ in $\Sigma$, and
$u_{n}\to u_{c}$ in $L^{p}$ for $p\in[2, 6)$. Combining with the weakly lower semi-continuous of norm, one has $u_{c}\in S(c)\cap B(r)$. Moreover, by Lemma
\ref{lem5} and the weakly lower semi-continuous of equivalent norm, we obtain that
\begin{align*}
E_{\Omega}(u_{c})\leq\liminf\limits_{n\to\infty} E_{\Omega}(u_{n})=m(c,r)\leq E_{\Omega}(u_{c}),
\end{align*}
which yields
$$
E_{\Omega}(u_{c})=m(c,r)
$$
and $u_{n}\to u_{c}$ in $\Sigma$. Thus, $u_{c}$ is a minimizer of $m(c,r)$ and any minimizing sequence of $m(c,r)$ is precompact.
Next, in order to claim that $u_{c}$ is a critical
point of $E_{\Omega}|_{S(c)}$. We need to show that $u_{c}$ is not on the boundary of $S(c)\cap B(r)$. Indeed, if $S(c)\cap(B(r)\backslash B(cr))=\emptyset$ with $c<1$,
it is easy to see that $\mathcal{M}_{c}^{r}\subseteq B(cr)$, so $u_{c}\notin S(c)\cap\partial B(r)$,
where $$
\partial B(r):=\{u\in \Sigma: |\nabla u|_2^2+|xu|_2^2=r\}.
$$
On the other hand, if $S(c)\cap(B(r)\backslash B(cr))\neq\emptyset$ with $c<1$, we claim
 that there exists some constant
$\widetilde{c}(r)>0$ such that for any $0<c<c_0:=\min\{1, \frac{r}{r_1}, \widetilde{c}(r)\}$,
\begin{align}\label{Q1}
\inf\limits_{u\in S(c)\cap B(\frac{C_{*}}{5C^{*}}cr)}E_{\Omega}(u)<\inf\limits_{u\in S(c)\cap(B(r)\setminus B(cr))}E_{\Omega}(u),
\end{align}
where $C_{*}$, $C^{*}$ are defined in Lemma \ref{lem5}.
In fact, by Lemmas \ref{lem2}, \ref{lem5}, and \eqref{GN}, one has
\begin{align}\label{J1}
E_{\Omega}(u)&=\frac{1}{2}\int_{\mathbb{R}^3}|\nabla u|^2dx+\int_{\mathbb{R}^3}V(x)|u|^2dx+\frac{1}{2}N(u)-\Omega\int_{\mathbb{R}^3}\overline{u} L_{z}udx\notag\\
&\geq \frac{C_{*}}{2}\int_{\mathbb{R}^3}|\nabla u|^2+|x|^2|u|^2dx-\frac{\Lambda}{2}\int_{\mathbb{R}^3}|u|^4dx\notag\\
&\geq \frac{C_{*}}{2}\int_{\mathbb{R}^3}|\nabla u|^2+|x|^2|u|^2dx-\frac{\Lambda}{2}\mathcal{C}_4^4
\left(\int_{\mathbb{R}^3}|\nabla u|^2+|x|^2|u|^2dx\right)^{\frac{3}{2}}c^{\frac{1}{2}}
\end{align}
for any $u\in S(c)$. Further, by Lemmas \ref{lem2}, \ref{lem5}, and Sobolev embedding inequality, we also have
\begin{align}\label{J2}
E_{\Omega}(u)&=\frac{1}{2}\int_{\mathbb{R}^3}|\nabla u|^2dx+\int_{\mathbb{R}^3}V(x)|u|^2dx+\frac{1}{2}N(u)-\Omega\int_{\mathbb{R}^3}\overline{u} L_{z}udx\notag\\
&\leq \frac{C^{*}}{2}\int_{\mathbb{R}^3}|\nabla u|^2+|x|^2|u|^2dx+\frac{\Lambda}{2}\int_{\mathbb{R}^3}|u|^4dx\notag\\
&\leq \frac{C^{*}}{2}\int_{\mathbb{R}^3}|\nabla u|^2+|x|^2|u|^2dx+\frac{\Lambda C}{2}\left(\int_{\mathbb{R}^3}|\nabla u|^2+|x|^2|u|^2dx\right)^2.
\end{align}
Let $$f_{c}(\eta):=\frac{C_{*}}{2}\eta-\frac{\Lambda}{2}\mathcal C_4^4c^{\frac{1}{2}}\eta^{\frac{3}{2}}$$ and $$g_{c}(\eta):=\frac{C^{*}}{2}\eta+\frac{\Lambda C}{2}\eta^2,$$
where $\eta:=\int_{\mathbb{R}^3}|\nabla u|^2+|x|^2|u|^2dx$ and $C_{*}$, $C^{*}$ are defined in Lemma \ref{lem5}. Thus, for any $\eta\in(0,r)$, we obtain that
$$
f_{c}(\eta)\leq E_{\Omega}(u)\leq g_{c}(\eta).
$$
By the definition of $g_{c}(\eta)$, we know that $g_{c}(\eta)$ is monotonically increasing about $\eta\in(0,r)$. Thus, it is sufficient to prove that there exists $\widetilde{c}(r)>0$ such that for any $0<c<\widetilde{c}(r)$,
\begin{align}\label{Q}
g_{c}\left(\frac{C_{*}}{5C^{*}}cr\right)<\inf\limits_{\eta\in(cr, r)}f_{c}(\eta).
\end{align}
Noticing that for any $\eta\in (cr,r)$, one has
\begin{align}\label{E>0}
f_{c}(\eta)&=\frac{C_{*}}{2}\eta-\frac{\Lambda}{2}\mathcal C_4^4c^{\frac{1}{2}}\eta^{\frac{3}{2}}\notag\\
&\geq\frac{C_{*}}{2}\eta\left(1-\frac{\Lambda}{C_{*}}\mathcal C_4^4c^{\frac{1}{2}}r^{\frac{1}{2}}\right)\notag\\
&>\frac{C_{*}}{4}\eta>\frac{C_{*}}{4}cr>0
\end{align}
provided that $c<c_1(r):=\frac{C_{*}^2}{4\Lambda^2\mathcal C_4^8r}$. We further obtain that
\begin{align*}
g_{c}\left(\frac{C_{*}}{5C^{*}}cr\right)=\frac{C_{*}cr}{10}+\frac{\Lambda C}{50}\left(\frac{C_{*}}{C^{*}}\right)^2c^2r^2<\frac{C_{*}cr}{8}
+\frac{C_{*}cr}{8}=\frac{C_{*}}{4}cr,
\end{align*}
if $c<c_2(r):=\frac{25(C^{*})^2}{4\Lambda CC_{*}r}$. Thus we conclude that \eqref{Q} holds for $c<\widetilde{c}(r)=\min\{c_1(r),c_2(r)\}$.
Consequently, choosing
\begin{align}\label{CS}
c_0:=\min\{1, \frac{r}{r_1}, \widetilde{c}(r)\},
\end{align}
we can infer that \eqref{Q1} holds for all $0<c<c_0$.
Now, in order to guarantee $u_c\notin S(c)\cap\partial B(r)$,
we also prove that $\mathcal{M}_{c}^{r}\subseteq B(cr)$ with $c<c_0$. Assume by contradiction that there is a $\varphi\in \mathcal{M}_{c}^{r}$ but
$\varphi\notin B(cr)$. It follows from \eqref{Q1} that for any $c<c_0$,
\begin{align*}
m(c,r)&\leq\inf\limits_{u\in S(c)\cap B(\frac{C_{*}}{5C^{*}}cr)}E_{\Omega}(u)\\
&<\inf\limits_{u\in S(c)\cap(B(r)\setminus B(cr))}E_{\Omega}(u)\\
&\leq E_{\Omega}(\varphi)=m(c,r),
\end{align*}
which yields a contradiction. Thus, we have $\mathcal{M}_{c}^{r}\subseteq B(cr)$, and $u_{c}$ is an interior local minimizer for $m(c,r)$. Namely, $u_{c}$ is a
critical point of $E_{\Omega}|_{S(c)}$ and there exists a corresponding Lagrange multiplier $\mu_{c}\in\mathbb R$.

In what follows, we estimate the bound of the multiplier $\mu_{c}$, that is, \eqref{mu-c-bound}.
We first give an upper bound of $m(c,r)$ from the spectral theory. Recall that $\Phi_c(x)=\sqrt{c}\Phi_0(x)\in S(c)$ and
 $\Phi_c\in S(c)\cap B(r)$ for $0<c<\frac{r}{r_1}$.
Moreover, one has $\int_{\mathbb{R}^3}\Omega\overline{\Phi_c}L_{z}\Phi_c dx=0$ as $\Phi_c$ is real valued. By Lemma \ref{lem2}, \eqref{GN} and \eqref{T}, we deduce that
\begin{align}\label{SJ}
m(c,r)\leq E_{\Omega}(\Phi_c)&=\frac{1}{2}\int_{\mathbb{R}^3}|\nabla \Phi_c|^2dx+\int_{\mathbb{R}^3}V(x)|\Phi_c|^2dx+N(\Phi_c)
-\Omega\int_{\mathbb{R}^3}\overline{\Phi_c} L_{z}\Phi_c dx\notag\\
&\leq c\left[\frac{1}{2}\int_{\mathbb{R}^3}|\nabla \Phi_0|^2dx+\int_{\mathbb{R}^3}V(x)|\Phi_0|^2dx\right]
+\frac{\Lambda}{2}\int_{\mathbb{R}^3}|\Phi_c|^4dx\notag\\
&\leq \zeta^0c+\sqrt{2}\Lambda \mathcal{C}_4^4\left(\frac{1}{2}\int_{\mathbb{R}^3}|\nabla \Phi_0|^2dx\right)^{\frac{3}{2}}c^2\notag\\
&\leq \zeta^0c+\sqrt{2}\Lambda \mathcal{C}_4^4(\zeta^0)^{\frac{3}{2}}c^2,
\end{align}
where $\zeta^0=\frac{1}{2}\sum\limits_{j=1}\limits^3\gamma_{j}$.
Note that $(u_{c}, \mu_{c})\in \mathcal{M}_{c}^{r}\times\mathbb{R}$ weakly solves problem \eqref{EQ}.
From Lemma \ref{lem2}, \eqref{GN} and \eqref{SJ}, we derive that
\begin{align}\label{Q2'}
\mu_{c}|u_{c}|_2^2&=-\frac{1}{2}\int_{\mathbb{R}^3}|\nabla u_{c}|^2dx-\int_{\mathbb{R}^3}V(x)|u_{c}|^2dx-N(u_{c})
+\Omega\int_{\mathbb{R}^3}\overline{u}_{c} L_{z}u_{c}dx\notag\\
&=-E_{\Omega}(u_{c})-\frac{1}{2}N(u_{c})\notag\\
&\geq-m(c,r)-\frac{\Lambda}{2}\int_{\mathbb{R}^3}|u_{c}|^4dx\notag\\
&\geq-m(c,r)-\frac{\Lambda}{2}\mathcal C_4^4\left(\int_{\mathbb{R}^3}|\nabla u_{c}|^2+|x|^2|u_{c}|^2dx\right)^{\frac{3}{2}}c^{\frac{1}{2}}\notag\\
&\geq-\zeta^0c-\sqrt{2}\Lambda \mathcal{C}_4^4(\zeta^0)^{\frac{3}{2}}c^2-\frac{\Lambda}{2}\mathcal C_4^4r^{\frac{3}{2}}c^2
\end{align}
due to $u_{c}\in\mathcal{M}_{c}^{r}\subseteq B(cr)$.

On the other hand, by Lemmas \ref{lem2}, \ref{lem5} and \eqref{GN}, we have
\begin{align*}
\mu_{c}|u_{c}|_2^2&=-\frac{1}{2}\int_{\mathbb{R}^3}|\nabla u_{c}|^2dx-\int_{\mathbb{R}^3}V(x)|u_{c}|^2dx-N(u_{c})
+\Omega\int_{\mathbb{R}^3}\overline{u}_{c} L_{z}u_{c}dx\\
&\leq -\frac{C_{*}}{2}\int_{\mathbb{R}^3}|\nabla u_{c}|^2+|x|^2|u_{c}|^2dx+\Lambda\int_{\mathbb{R}^3}|u_{c}|^4dx\\
&\leq -\frac{C_{*}}{2}\int_{\mathbb{R}^3}|\nabla u_{c}|^2+|x|^2|u_{c}|^2dx
+\Lambda \mathcal C_4^4\left(\int_{\mathbb{R}^3}|\nabla u_{c}|^2+|x|^2|u_{c}|^2dx\right)^{\frac{3}{2}}c^{\frac{1}{2}}\\
&=\int_{\mathbb{R}^3}|\nabla u_{c}|^2+|x|^2|u_{c}|^2dx\left[-\frac{C_{*}}{2}+\Lambda \mathcal C_4^4\left(\int_{\mathbb{R}^3}|\nabla u_{c}|^2+|x|^2|u_{c}|^2dx\right)^{\frac{1}{2}}c^{\frac{1}{2}}\right]\\
&\leq\left[-\frac{C_{*}}{2}+\Lambda \mathcal C_4^4r^{\frac{1}{2}}c^{\frac{1}{2}}\right]\int_{\mathbb{R}^3}|\nabla u_{c}|^2+|x|^2|u_{c}|^2dx,
\end{align*}
where $C_{*}=\min\{\frac{(\Omega^{*})^2-\Omega^2}{(\Omega^{*})^2+\Omega^2}, \frac{(\Omega^{*})^2-\Omega^2}{2}, \gamma_3^2\}$ defined in Lemma \ref{lem5}.
Noticing that $0<c<c_0$ (see \eqref{CS}), we obtain that $-\frac{C_{*}}{2}+\Lambda \mathcal C_4^4r^{\frac{1}{2}}c^{\frac{1}{2}}<0$. Then we obtain from \eqref{T} that
\begin{align*}
\mu_{c}|u_{c}|_2^2\leq3\left[-\frac{C_{*}}{2}+\Lambda \mathcal C_4^4r^{\frac{1}{2}}c^{\frac{1}{2}}\right]c<0
\end{align*}
for $0<c<c_0$.
Combining this with \eqref{Q2'}, we conclude that for any $0<c<c_0$,
\begin{align*}
-\zeta^0-\delta(r)c\leq\mu_{c}\leq3\left[-\frac{C_{*}}{2}+\Lambda \mathcal C_4^4r^{\frac{1}{2}}c^{\frac{1}{2}}\right],
\end{align*}
where $\delta(r):=\sqrt{2}\Lambda \mathcal{C}_4^4(\zeta^0)^{\frac{3}{2}}+\frac{\Lambda}{2}\mathcal C_4^4r^{\frac{3}{2}}$. This completes the proof.\qed
\\~

\subsection{Existence of a second critical point of mountain pass type}

Motivated by \cite{Jeanjean1997,Luo2021}, in this subsection, we prove that the existence of a mountain pass solution to \eqref{EQ}.
The proof  is divided into several steps:

\textbf{\emph {Step} 1}: $E_{\Omega}(u)$ has a local mountain pass geometry on $S(c)$.

By Theorem \ref{thm1}, there exists $u_{c}\in\mathcal{M}_{c}^{r}$ such that $E(u_{c})=m(c,r)$. We consider the scaling given by
\begin{align*}
u_{t}(x)=t^{\frac{5}{4}}u_{c}(tx_1, tx_2, t^{\frac{1}{2}}x_3),\,\,\,\,t>0,
\end{align*}
which implies that $u_{t}\in S(c)$ for all $t>0$.
Then, by a change of variable and $\widehat{u_{t}^2}=\mathcal{F}(u_{t})^2(\xi)=\mathcal{F}u^2(\frac{\xi_1}{t}, \frac{\xi_2}{t}, \frac{\xi_3}{t})$, the energy functional
rescales as
\begin{align*}
E_{\Omega}(u_{t})=&\frac{t^2}{2}\int_{\mathbb{R}^3}|\nabla_{x_1, x_2}u_{c}|^2dx+\frac{t}{2}\int_{\mathbb{R}^3}|\nabla_{x_3}u_{c}|^2dx
+\frac{1}{2}\int_{\mathbb{R}^3}\left(\gamma_1^2\frac{x_1^2}{t^2}+\gamma_2^2\frac{x_2^2}{t^2}+\gamma_3^2\frac{x_3^2}{t}\right)|u_{c}|^2dx\\
&+\frac{t^{\frac{5}{2}}}{2}\frac{1}{(2\pi)^3}\int_{\mathbb{R}^3}\left(\lambda_1
+\frac{4\pi}{3}\lambda_2\frac{2t\xi_3^2-t^2\xi_1^2-t^2\xi_2^2}{t^2\xi_1^2+t^2\xi_2^2+t\xi_3^2}\right)|\widehat{u_{c}^2}|^2d\xi
-\int_{\mathbb{R}^3}\Omega\overline{u_{c}}L_{z}u_{c}dx.
\end{align*}
 Recall that $(\lambda_1, \lambda_2)\in D_{ur}$, then if $\lambda_2>0$, we have
\begin{align*}
\lim\limits_{t\to\infty}\lambda_1+\frac{4\pi}{3}\lambda_2\frac{2t\xi_3^2-t^2\xi_1^2-t^2\xi_2^2}{t^2\xi_1^2+t^2\xi_2^2+t\xi_3^2}=\lambda_1-\frac{4\pi}{3}\lambda_2<0.
\end{align*}
It yields that $\lim\limits_{t\to\infty}E_{\Omega}(u_{t})=-\infty$ thanks to the Lebesgue's theorem.
On the other hand, when $\lambda_2<0$, the same conclusion is obtained by choosing the scaling
\begin{align*}
\widetilde{u}_{t}(x)=t^{\frac{5}{4}}u_{c}(t^{\frac{3}{4}}x_1, t^{\frac{3}{4}}x_2, tx_3),\,\,\,\,t>0.
\end{align*}
Thus, in both cases, the class of paths $\Gamma(c)$ can be defined as follows
\begin{align}\label{O3}
\Gamma(c):=\{g\in C([0,1], S(c)): g(0)=u_{c},\,\, E_{\Omega}(g(1))<0\}.
\end{align}
Note that $g(s)=(1+(t-1)s)^{\frac{5}{4}}u_{c}((1+(t-1)s)x_1, (1+(t-1)s)x_2, (1+(t-1)s^{\frac{1}{2}})x_3)\in\Gamma(c)$ with $t\gg1$ if $\lambda_2>0$
and the case of $\lambda_2<0$ is similar. Then $\Gamma(c)\neq\emptyset$
and we can
introduce a min-max value
\begin{align}\label{O2}
\gamma(c):=\inf\limits_{g\in\Gamma(c)}\max\limits_{0\leq s\leq1}E_{\Omega}(g(s)).
\end{align}
Next, by \eqref{E>0} and \eqref{O3}-\eqref{O2}, we deduce that
\begin{align}\label{O1}
\gamma(c)>\max\{E_{\Omega}(u_{c}), E_{\Omega}(g(1))\}>0.
\end{align}

\textbf{\emph{Step} 2}: Existence of bounded special Palais-Smale sequences.

In order to obtain a special Palais-Smale sequence, we introduce the auxiliary functional
\begin{align*}
\widetilde{E}_{\Omega}: S(c)\times\mathbb{R}\rightarrow\mathbb{R}\,\,\,\,(u,t)\rightarrow E_{\Omega}(\kappa(u,t)),
\end{align*}
where $\kappa(u,t):=e^{\frac{3}{2}t}u(e^{t}x)$. More precisely, we have
\begin{align*}
\widetilde{E}_{\Omega}(u,t)=E_{\Omega}(\kappa(u,t))&=E_{\Omega}(e^{\frac{3}{2}t}u(e^{t}x))\\
&=\frac{e^{2t}}{2}\int_{\mathbb{R}^3}|\nabla u|^2dx+\frac{1}{e^{2t}}\int_{\mathbb{R}^3}V(x)|u|^2dx
+\frac{e^{3t}}{2}N(u)-\Omega\int_{\mathbb{R}^3}\overline{u} L_{z}udx.
\end{align*}
Define the set of paths
\begin{align*}
\widetilde{\Gamma}(c):=\{\widetilde{g}\in C([0,1]\to S(c)\times\mathbb{R}): \widetilde{g}(0)=(u_{c}, 0), \widetilde{g}(1)=(g(1), 0)\},
\end{align*}
and the corresponding of minimax value
\begin{align*}
\widetilde{\gamma}(c):=\inf\limits_{\widetilde{g}\in\widetilde{\Gamma}(c)}\max\limits_{s\in[0,1]}\widetilde{E}_{\Omega}(\widetilde{g}(s)).
\end{align*}
We can infer that $\widetilde{\gamma}(c)=\gamma(c)$. Indeed, since $\Gamma(c)\subseteq\widetilde{\Gamma}(c)$, it directly yields that
$\widetilde{\gamma}(c)\leq\gamma(c)$. Then 
we only need to prove that $\widetilde{\gamma}(c)\geq\gamma(c)$. This follows from the claim: for any $\widetilde{g}\in\widetilde{\Gamma}(c)$,
there is a $g\in\Gamma(c)$ such that
\begin{align}\label{Q3}
\max\limits_{s\in[0,1]}\widetilde{E}_{\Omega}(\widetilde{g}(s))=\max\limits_{s\in[0,1]}E_{\Omega}(g(s)).
\end{align}
Setting $\widetilde{g}(s):=(g_{1}(s), t)\in S(c)\times\mathbb{R}$, we get
\begin{align*}
\widetilde{E}_{\Omega}(\widetilde{g}(s))=\widetilde{E}_{\Omega}(g_{1}(s),t)=E_{\Omega}(\kappa(g_{1}(s),t))
\end{align*}
for all $s\in[0,1]$, which yields that $g(s):=\kappa(g_1(s), t)\in \Gamma(c)$. Namely, \eqref{Q3} holds and then $\widetilde{\gamma}(c)=\gamma(c)$.

Next, we give the following useful Lemma, which is crucial for
 obtaining the existence of bounded Palais-Smale sequence of $E_{\Omega}$ restricted to $S(c)$ at the level $\gamma(c)$,
In particular, it can be obtained by \cite[Lemma 2.3]{Jeanjean1997}.

Denote $X:=\Sigma\times\mathbb{R}$ endowed
with the norm $\|\cdot\|_{X}^2=\|\cdot\|_{\Sigma}^2+|\cdot|_{\mathbb{R}}^2$ where $|r|_{\mathbb{R}}=|r|$ for $r\in\mathbb{R}$, and $X^{-1}$ the dual space of $X$.
\begin{lem}
\label{lem7}
Suppose that $(\lambda_1, \lambda_2)\in D_{ur}$ and let $\varepsilon>0$ and $\widetilde{g}_0\in \widetilde{\Gamma}(c)$ be such that
\begin{align*}
\max\limits_{s\in[0,1]}\widetilde{E}_{\Omega}(\widetilde{g}_0(s))\leq\widetilde{\gamma}(c)+\varepsilon.
\end{align*}
Then there exists a pair of $(u_0,t_0)\in S(c)\times\mathbb{R}$ such that
\begin{itemize}
  \item  $\widetilde{E}_{\Omega}(u_0,t_0)\in [\widetilde{\gamma}(c)-\varepsilon, \widetilde{\gamma}(c)+\varepsilon]$;

  \item  $\min\limits_{s\in[0,1]}\|(u_0,t_0)-\widetilde{g}_0(s)\|_{X}\leq\sqrt{\varepsilon}$;

  \item  $\|(\widetilde{E}_{\Omega}\mid_{S(c)\times\mathbb{R}})'(u_0,t_0)\|_{X^{-1}}\leq2\sqrt{\varepsilon}$, that is
  $$
  |\langle\widetilde{E}'_{\Omega}(u_0, t_0), z\rangle_{X^{-1}\times X}|\leq2\sqrt{\varepsilon}\|z\|_{X}
  $$
  for all $z\in \widetilde{T}_{(u_0,t_0)}:=\{(z_1,z_2)\in X: \langle u_0, z_1\rangle=0\}$.
\end{itemize}
\end{lem}

Based on the above result, we can obtain the special Palais-Smale sequence for $\gamma(c)$.
\begin{lem}\label{lem8}
Suppose that $(\lambda_1, \lambda_2)\in D_{ur}$, $0<\Omega<\Omega^{*}$, and $0<c<c_0$ for $c_0$ defined in Theorem {\rm\ref{thm1}}.
Then there exists a sequence $\{v_{n}\}\subseteq S(c)$ such that
\begin{align}\label{Q4}
\begin{cases}
E_{\Omega}(v_{n})\to \gamma(c),\\
(E_{\Omega}\mid_{S(c)})'(v_{n})\to0,\\
Q(v_{n})\to0,
\end{cases}
\end{align}
as $n\to\infty$.
\end{lem}
\begin{prf}
By the definition of $\gamma(c)$, there exists a $g_{n}\in \Gamma(c)$ such that
$$
\max\limits_{s\in[0,1]}E_{\Omega}(g_{n}(s))\leq\gamma(c)+\frac{1}{n}.
$$
Note that $\widetilde{\gamma}(c)=\gamma(c)$. Letting $\widetilde{g}_{n}(s)=(g_{n}(s), 0)\in \widetilde{\Gamma}(c)$,
one has $\max\limits_{s\in[0,1]}\widetilde{E}_{\Omega}(\widetilde{g}_{n}(s))\leq\widetilde{\gamma}(c)+\frac{1}{n}$. Then, by applying Lemma \ref{lem7},
we deduce that there exists a sequence $\{(u_{n}, t_{n})\}\subseteq S(c)\times\mathbb{R}$ such that
\begin{enumerate}[\indent$(i)$]
   \item
   $\widetilde{E}_{\Omega}(u_{n},t_{n})\in [\widetilde{\gamma}(c)-\frac{1}{n}, \widetilde{\gamma}(c)+\frac{1}{n}]$;\label{i}

   \item
    $\min\limits_{s\in[0,1]}\|(u_{n},t_{n})-\widetilde{g}_{n}(s)\|_{X}\leq\sqrt{\frac{1}{n}}$;\label{ii}

   \item
    $\|(\widetilde{E}_{\Omega}\mid_{S(c)\times\mathbb{R}})'(u_{n},t_{n})\|_{X^{-1}}\leq2\sqrt{\frac{1}{n}}$, that is
  $$
  |\langle\widetilde{E}'_{\Omega}(u_{n}, t_{n}), z\rangle_{X^{-1}\times X}|\leq2\sqrt{\frac{1}{n}}\|z\|_{X}
  $$
  for all $z\in \widetilde{T}_{(u_{n},t_{n})}:=\{(z_1,z_2)\in X: \langle u_{n}, z_1\rangle=0\}$.\label{iii}
\end{enumerate}
Denote $v_{n}=\kappa(u_{n}, t_{n})$, then we claim that $\{v_{n}\}\subseteq S(c)$ satisfies \eqref{Q4}.
In fact, since $E_{\Omega}(v_{n})=E_{\Omega}(\kappa(u_{n},t_{n}))=\widetilde{E}_{\Omega}(u_{n}, t_{n})$, it follows from $(\ref{i})$ that $E_{\Omega}(v_{n})\to \gamma(c)$
as $n\to\infty$. By direct calculations, we have
$$
\frac{\partial\widetilde{E}_{\Omega}(u_{n},t_{n})}{\partial t_{n}}=2Q(v_{n}).
$$
Let $\frac{\partial\widetilde{E}_{\Omega}(u_{n},t_{n})}{\partial t_{n}}=\langle\widetilde{E}'_{\Omega}(u_{n},t_{n}), (0,1)\rangle_{X^{-1}\times X}$.
Since $(0,1)\in \widetilde{T}_{(u_{n},t_{n})}$, then $(\ref{iii})$ implies that $Q(v_{n})\to 0$ as $n\to \infty$.
Moreover, in order to verify
that $(E_{\Omega}\mid_{S(c)})'(v_{n})\to0$ as $n\to\infty$, it suffices to claim that for $n\in\mathbb{N}$ sufficiently large that
\begin{align}\label{Q5}
|\langle E'_{\Omega}(v_{n}), \varphi\rangle|\leq\frac{2\sqrt{2}}{\sqrt{n}}\|\varphi\|_{\Sigma}\,\,\,{\rm{for}}\,\,{\rm{all}}\,\,\,\varphi\in T_{v_{n}},
\end{align}
where $T_{v_{n}}:=\{\varphi\in\Sigma: \langle v_{n}, \varphi\rangle=0\}$.
Now we prove the above claim. Set $\widetilde{\varphi}=\kappa(\varphi, -t_{n})$ for any $\varphi\in T_{v_{n}}=\{\varphi\in\Sigma: \langle v_{n}, \varphi\rangle=0\}$, then one has
\begin{align*}
\langle E'_{\Omega}(v_{n}), \varphi\rangle=\langle\widetilde{E}_{\Omega}'(u_{n}, t_{n}), (\widetilde{\varphi}, 0)\rangle.
\end{align*}
Since $\int_{\mathbb{R}^3}u_{n}\widetilde{\varphi}dx=\int_{\mathbb{R}^3}v_{n}\varphi dx$,
we obtain that $(\widetilde{\varphi}, 0)\in \widetilde{T}_{(u_{n}, t_{n})}\Leftrightarrow \varphi\in T_{v_{n}}$. It follows from $(\ref{iii})$ that
\begin{align*}
|\langle E_{\Omega}'(v_{n}), \varphi\rangle|
=|\langle\widetilde{E}_{\Omega}'(u_{n}, t_{n}), (\widetilde{\varphi}, 0)\rangle|\leq2\sqrt{\frac{1}{n}}\|(\widetilde{\varphi}, 0)\|_{X}.
\end{align*}
Consequently, it is equivalent to prove that
$$
\|(\widetilde{\varphi}, 0)\|_{X}\leq\sqrt{2}\|\varphi\|_{\Sigma}
$$
for $n\in\mathbb{N}$ sufficiently large. From $(\ref{ii})$, one has
\begin{align*}
|t_{n}|=|t_{n}-0|\leq\min\limits_{s\in[0,1]}\|(u_{n}, t_{n})-(g_{n}(s), 0)\|_{X}\leq\sqrt{\frac{1}{n}},
\end{align*}
which yields that
\begin{align*}
\|(\widetilde{\varphi}, 0)\|_{X}^2=\|\widetilde{\varphi}\|_{\Sigma}^2
&=\int_{\mathbb{R}^3}|\nabla\widetilde{\varphi}|^2+|x|^2|\widetilde{\varphi}|^2+|\widetilde{\varphi}|^2dx\\
&=e^{-2t_{n}}\int_{\mathbb{R}^3}|\nabla \varphi|^2dx+e^{2t_{n}}\int_{\mathbb{R}^3}|x|^2|\varphi|^2dx+\int_{\mathbb{R}^3}|\varphi|^2dx\\
&\leq2\|\varphi\|_{\Sigma}^2
\end{align*}
for $n\in\mathbb{N}$ sufficiently large. This completes the proof.
\end{prf}

\textbf{{\emph {Step}} 3}: Existence of a mountain pass solution.

Now, we give the compactness of the Palais-Smale sequence obtained in Lemma \ref{lem8}.
\begin{lem}\label{lem9}
Suppose that $(\lambda_1, \lambda_2)\in D_{ur}$, $0<\Omega<\frac{\sqrt{5}}{3}\Omega^{*}$, and $0<c<c_0$ for $c_0$ defined in Theorem {\rm\ref{thm1}}.
Let $\{v_{n}\}\subseteq S(c)$ obtained in Lemma {\rm\ref{lem8}}. Then, up to a subsequence,
\begin{enumerate}[\indent$(i)$]
\item
   there is $v\in\Sigma$ such that $v_{n}\rightharpoonup v$ weakly in $\Sigma$;
\item
   there is $\{\mu_{n}\}\subseteq\mathbb{R}$ such that ${\rm Re}\left[E_{\Omega}'(v_{n})+\mu_{n}v_{n}\right]\to 0$ in $\Sigma^{-1}$ as $n\to \infty$;
\item
   there is a constant $\widehat{\mu}\in\mathbb{R}$ such that $\mu_{n}\to\widehat{\mu}$ in $\mathbb{R}$ and ${\rm Re}\left[E_{\Omega}'(v)+\widehat{\mu}v\right]= 0$ in $\Sigma^{-1}$ as $n\to \infty$.
\end{enumerate}
\end{lem}
\begin{prf}
To claim that $(i)$ holds, we first prove that $\{v_{n}\}$ is bounded in $\Sigma$.
Note that $Q(v_{n})\to 0$ as $n\to\infty$ and $0<\Omega<\frac{\sqrt{5}}{3}\Omega^{*}$, then it follows from \eqref{E1} that
\begin{align*}
&\ \ \gamma(c)+o_{n}(1)\\
&=E_{\Omega}(v_n)-\frac{2}{3}Q(v_n)+o_n(1)\\
&=\frac{1}{6}\int_{\mathbb{R}^3}|\nabla v_n|^2dx+\frac{5}{6}\int_{\mathbb{R}^3}
\sum\limits_{j=1}\limits^{3}\gamma_{j}^2x_{j}^2|v_n|^2dx
-\int_{\mathbb{R}^3}\Omega\overline{v_n} L_{z}v_ndx+o_n(1)\\
&\geq(\frac{1}{6}-\varepsilon)\int_{\mathbb{R}^3}|\nabla v_n|^2dx
+\int_{\mathbb{R}^3}\left[(\frac{5}{6}\gamma_1^2-
\frac{\Omega^2}{4\varepsilon})x_1^2
+(\frac{5}{6}\gamma_2^2-\frac{\Omega^2}{4\varepsilon})
x_2^2+\frac{5}{6}\gamma_3^2x_3^2\right]|v_n|^2dx+o_n(1).
\end{align*}
In view of the proof Lemma \ref{lem5}, we can choose $\varepsilon=\varepsilon_1\in(\frac{3\Omega^2}{10(\Omega^{*})^2}, \frac{1}{6})$ as $0<\Omega<\frac{\sqrt{5}}{3}\Omega^{*}$ and then
\begin{align*}
\widetilde{C}_1:=\frac{1}{6}-\varepsilon_1>0, \,\,\,\widetilde{C}_2:=\frac{5}{6}(\Omega^{*})^2-
\frac{\Omega^2}{4\varepsilon_1}>0.
\end{align*}
Letting
\begin{align*}
\widetilde{C}:=\min\{\widetilde{C}_1, \widetilde{C}_2, \frac{5}{6}\gamma_3^2\},
\end{align*}
we can obtain that
\begin{align}\label{bounded-PS-sequence}
\gamma(c)+o_{n}(1)\geq\widetilde{C}\int_{\mathbb{R}^3}|\nabla v_n|^2+|x|^2|v_n|^2dx.
\end{align}
Thus, $\{v_{n}\}$ is bounded in $\Sigma$. Then, up to a subsequence, $(i)$ holds. Next, to obtain $(ii)$, it suffices to
claim that there exists $C>0$ independent of $n\in\mathbb{N}$ such that
\begin{align}\label{Q6}
|{\rm{Re}}\langle E_{\Omega}'(v_{n})-\frac{1}{c}\langle E_{\Omega}'(v_{n}), v_{n}\rangle v_{n}, z\rangle_{\Sigma^{-1}\times\Sigma}|\leq\frac{C}{\sqrt{n}}
\|z\|_{\Sigma}
\end{align}
for any $z\in \Sigma$. More specifically, we know that
\begin{align}\label{Q6'}
{\rm{Re}}\langle E_{\Omega}'(v_{n})-\frac{1}{c}\langle E_{\Omega}'(v_{n}), v_{n}\rangle v_{n}, z\rangle_{\Sigma^{-1}\times\Sigma}={\rm{Re}}\langle E_{\Omega}'(v_{n}), \widetilde{z}_{n}\rangle_{\Sigma^{-1}\times\Sigma}
\end{align}
for any $z\in \Sigma$, where $\widetilde{z}_{n}:=z-\frac{1}{c}\langle v_{n},z\rangle_{L^2(\mathbb{R}^3)} v_{n}$. It is easy to see that $\widetilde{z}_{n}\in T_{v_{n}}$.
Since $(E_{\Omega}\mid_{S(c)})'(v_{n})\to 0$ as $n\to\infty$, it follows from \eqref{Q6'} that
\begin{align*}
|{\rm{Re}}\langle E_{\Omega}'(v_{n})-\frac{1}{c}\langle E_{\Omega}'(v_{n}), v_{n}\rangle v_{n}, z\rangle_{\Sigma^{-1}\times\Sigma}|\leq\frac{2\sqrt{2}}{\sqrt{n}}\|\widetilde{z}_{n}\|_{\Sigma}.
\end{align*}
By simple calculations, one has
\begin{align*}
\|\widetilde{z}_{n}\|_{\Sigma}&=\|z-\frac{1}{c}\langle v_{n}, z\rangle_{L^2(\mathbb{R}^3)} v_{n}\|_{\Sigma}\\
&\leq\|z\|_{\Sigma}+\frac{1}{c}|\langle v_{n}, z\rangle_{L^2(\mathbb{R}^3)}|\|v_{n}\|_{\Sigma}\\
&\leq\|z\|_{\Sigma}+\frac{|v_{n}|_2}{c}|z|_2|\|v_{n}\|_{\Sigma}\\
&\leq\|z\|_{\Sigma}\left(1+\frac{|v_{n}|_2}{c}\|v_{n}\|_{\Sigma}\right).
\end{align*}
By applying the fact that $|v_{n}|_2^2=c$ and $\{v_{n}\}$ is bounded in $\Sigma$, we deduce that \eqref{Q6} holds. Namely,
\begin{align*}
&{\rm{Re}}\left[\frac{1}{2}\int_{\mathbb{R}^3}\nabla v_{n}\cdot\nabla\overline{\varphi}+\int_{\mathbb{R}^3}V(x)v_{n}\overline{\varphi}dx
+\frac{1}{2}\int_{\mathbb{R}^3}\lambda_1|v_{n}|^2v_{n}\overline{\varphi}+\lambda_2(K\ast|v_{n}|^2)v_{n}\overline{\varphi}dx\right.\\
&\left.\,\,\,\qquad\,\,-\int_{\mathbb{R}^3}\Omega\overline{\varphi} L_{z}v_{n}dx+\mu_{n}\int_{\mathbb{R}^3}v_{n}\overline{\varphi}
dx\right]\to0
\end{align*}
as $n\to\infty$, where
\begin{align*}
\mu_{n}&=-\frac{1}{c}\langle E_{\Omega}'(v_{n}), v_{n}\rangle\\
&=-\frac{1}{c}\left[\frac{1}{2}\int_{\mathbb{R}^3}|\nabla v_{n}|^2dx+\int_{\mathbb{R}^3}V(x)|v_{n}|^2dx
+\frac{1}{2}N(v_{n})-\int_{\mathbb{R}^3}\Omega\overline{v_{n}} L_{z}v_{n}dx\right].
\end{align*}
 Thus, $(ii)$ holds. To the end, we prove that there is $\widehat{\mu}\in\mathbb{R}$
such that, up to a subsequence, $\mu_{n}\to\widehat{\mu}$ in $\mathbb{R}$. By Lemmas \ref{lem2} and \ref{lem5}, we deduce that
\begin{align*}
\mu_{n}c=-\left(\frac{1}{2}\int_{\mathbb{R}^3}|\nabla v_{n}|^2dx+\int_{\mathbb{R}^3}V(x)|v_{n}|^2dx
+\frac{1}{2}N(v_{n})-\int_{\mathbb{R}^3}\Omega\overline{v_{n}} L_{z}v_{n}dx\right),
\end{align*}
and then there is $\widehat{\mu}\in\mathbb{R}$ such that, up to a subsequence,
$\mu_{n}\to\widehat{\mu}$ in $\mathbb{R}$. Combining this with $(ii)$, it yields that ${\rm Re}\left[E_{\Omega}'(v)+\widehat{\mu}v\right]= 0$ in $\Sigma^{-1}$ as $n\to \infty$. This completes the proof.
\end{prf}

\noindent\textbf{Proof of Theorem \ref{thm1} (existence of a second critical point)}.
By the proof of Lemmas \ref{lem8}-\ref{lem9}, we obtain that a specific sequence $\{v_{n}\}\subseteq S(c)$ and $v_{n}\rightharpoonup v$ in $\Sigma$.
Then, by Lemma \ref{lem4}, we have that
\begin{align}\label{Q7}
v_{n}\to v\,\,\,{\rm{in}}\,\,\, L^{p}(\mathbb{R}^3), \,\,2\leq p<6.
\end{align}
Further, by Lemma \ref{lem9}-$(ii)$ and $(iii)$, we deduce that
\begin{align}\label{Q8}
{\rm{Re}}\langle E'(v_{n})+\widehat{\mu}v_{n}, v_{n}-v\rangle=o_{n}(1),\,\,\,
{\rm{Re}}\langle E'(v)+\widehat{\mu}v, v_{n}-v\rangle=0.
\end{align}
By \cite[Lemma 8.1]{Minimax} and \eqref{Q7} as well as Lemma \ref{lem5}, we can deduce from \eqref{Q8} that
\begin{align*}
o_{n}(1)&=\frac{1}{2}\int_{\mathbb{R}^3}|\nabla (v_{n}-v)|^2dx+\int_{\mathbb{R}^3}V(x)|v_{n}-v|^2dx-\int_{\mathbb{R}^3}\Omega(\overline{v_{n}-v}) L_{z}(v_{n}-v)dx\\
&\geq\frac{C_{*}}{2}\int_{\mathbb{R}^3}|\nabla(v_{n}-v)|^2+|x|^2|v_{n}-v|^2dx,
\end{align*}
which implies that $v_{n}\to v$ in $\Sigma$. Therefore, we obtain a mountain pass solution
$(v, \widehat{\mu})\in \Sigma\times\mathbb{R}$ to problem \eqref{EQ} with $|v|_2^2=c$. Taking into account \eqref{O1}, one has
\begin{align*}
E_{\Omega}(v)=\gamma(c)>E(u_{c})=m(c,r).
\end{align*}
This completes the proof.\qed

\subsection{Asymptotic properties and ground states}
In this subsection, we give proof of Theorems \ref{thm1''} and \ref{thm1'}.

\noindent\textbf{Proof of Theorem \ref{thm1''}}.
From Theorem \ref{thm1}, we know that for any fixed $r>0$ and every $0<c<c_0$,
$$
\mathcal{M}_{c}^{r}=\{u\in S(c)\cap B(r): E_{\Omega}(u)=m(c,r)\}\neq\emptyset.
$$
Letting $u\in\mathcal{M}_{c}^{r}$, we rewrite $u=u_1+iu_2$, where $u_1$ is the real part and $u_2$ is the imaginary part of $u$.
Motivated by \cite[Theorem 2]{Bellazzini2017}, we introduce the pure point spectrum $\lambda_{j}$ of the harmonic oscillator $-\frac{1}{2}\Delta+\frac{|x|^2}{2}$
and the corresponding eigenfunction given by Hermite function $\psi_{j}$ such that
\begin{align*}
-\frac{1}{2}\Delta\psi_{j}+\frac{|x|^2}{2}\psi_{j}&=\lambda_{j}\psi_{j},\,\,\,\int_{\mathbb{R}^3}|\psi_{j}|^2dx=1,\\
\lambda_{j}&\leq\lambda_{j+1},\,\,\,j=0,1,2,\cdots.
\end{align*}
As is known to all, $\{\psi_{j}\}$ is a Hilbert basis for $L^2(\mathbb{R}^3, \mathbb{R})$ and $\lambda_{j}=\frac{3}{2}+j$ (see \cite{Antonelli2015}).
By applying Fourier decomposition, one has
\begin{align}\label{SD}
u=\sum\limits_{j\geq0}\left(\int_{\mathbb{R}^3}u_1\psi_{j}dx\right)\psi_{j}+i\sum\limits_{j\geq0}\left(\int_{\mathbb{R}^3}u_2\psi_{j}dx\right)\psi_{j}
=\sum\limits_{j\geq0}\varrho_{j}\psi_{j}
\end{align}
with $\varrho_{j}=\int_{\mathbb{R}^3}u\psi_{j}dx$.
So
\begin{align}\label{DD}
c=|u|_2^2=\sum\limits_{j\geq0}\varrho_{j}\overline{\varrho}_{j}\int_{\mathbb{R}^3}|\psi_{j}|^2dx=\sum\limits_{j\geq0}|\varrho_{j}|^2,
\end{align}
where $\overline{\varrho}_{j}$ is the conjugate of $\varrho_{j}$.
In view of \eqref{SJ}, we obtain that
\begin{align}\label{Sup}
E_{\Omega}(u)=m(c,r)\leq \zeta^0c+\sqrt{2}\Lambda \mathcal{C}_4^4(\zeta^0)^{\frac{3}{2}}c^2,
\end{align}
where $\zeta^0=\frac{1}{2}\sum\limits_{j=1}\limits^3\gamma_{j}$.
Since $u\in\mathcal{M}_{c}^{r}\subseteq B(cr)$, it follows from \eqref{J1} that
\begin{align*}
\zeta^0c+\sqrt{2}\Lambda \mathcal{C}_4^4(\zeta^0)^{\frac{3}{2}}c^2&\geq E_{\Omega}(u)\\
&\geq \frac{C_{*}}{2}\int_{\mathbb{R}^3}|\nabla u|^2+|x|^2|u|^2dx-\frac{\Lambda}{2}\mathcal C_4^4c^2r^{\frac{3}{2}}\\
&=C_{*}\sum\limits_{j\geq0}\lambda_{j}|\varrho_{j}|^2-\frac{\Lambda}{2}\mathcal C_4^4r^{\frac{3}{2}}c^2,
\end{align*}
which yields that
\begin{align}\label{DD'}
\sum\limits_{j\geq0}\lambda_{j}
|\varrho_{j}|^2\leq\frac{1}{C_{*}}\left(\zeta^0c+\delta(r)c^2
\right).
\end{align}
Here $\delta(r):=\sqrt{2}\Lambda \mathcal{C}_4^4(\zeta^0)^{\frac{3}{2}}+\frac{\Lambda}{2}\mathcal C_4^4r^{\frac{3}{2}}$ and
$C_{*}=\min\{\frac{(\Omega^{*})^2-\Omega^2}{(\Omega^{*})^2+\Omega^2}, \frac{(\Omega^{*})^2-\Omega^2}{2}, \gamma_3^2\}$ defined in Lemma \ref{lem5}.
By the definition of $\|\cdot\|_{\dot{\Sigma}}$ and \eqref{SD}, \eqref{DD'}, we conclude that
\begin{align*}
\frac{1}{2}\|u-\varrho_0\psi_0\|^2_{\dot{\Sigma}}&=\sum\limits_{j\geq1}\lambda_{j}|\varrho_{j}|^2
\leq\frac{1}{C_{*}}\left(\zeta^0c+\delta(r)c^2\right),
\end{align*}
which implies that
$$
\sup\limits_{u\in \mathcal{M}_{c}^{r}}\|u-\varrho_0\psi_0\|^2_{\dot{\Sigma}}\leq\frac{2}{C_{*}}\left[\zeta^0c+\left(\sqrt{2}\Lambda \mathcal{C}_4^4(\zeta^0)^{\frac{3}{2}}+\frac{\Lambda}{2}\mathcal C_4^4r^{\frac{3}{2}}\right)c^2\right].
$$
Thus
$$\sup\limits_{u\in \mathcal{M}_{c}^{r}}\|u-\varrho_0\psi_0\|^2_{\dot{\Sigma}}=O(c+c^2).
$$
Further, we deduce from \eqref{DD} and \eqref{DD'} that
\begin{align*}
\inf\limits_{j\geq1}\lambda_{j}\left(1-\frac{|\varrho_0|^2}{c}\right)\leq\lambda_0\left(\frac{\zeta^0}{C_{*}\lambda_0}-\frac{|\varrho_0|^2}{c}\right)+\frac{1}{C_{*}}\delta(r)c.
\end{align*}
Setting $\zeta^0= C_{*}\lambda_0$, i.e., $\zeta^0= \frac{3C_{*}}{2}$, we have
\begin{align}
\inf\limits_{j\geq1}(\lambda_{j}-\lambda_0)\left(\frac{\zeta^0}{C_{*}\lambda_0}-\frac{|\varrho_0|^2}{c}\right)\leq\frac{1}{C_{*}}\delta(r)c.
\end{align}
Using \eqref{DD'} again, we conclude that
\begin{align*}
\frac{\frac{1}{2}\|u-\varrho_0\psi_0\|^2_{\dot{\Sigma}}}{c}&=\sum\limits_{j\geq1}\lambda_{j}\frac{|\varrho_{j}|^2}{c}\\
&=\sum\limits_{j\geq0}\lambda_{j}\frac{|\varrho_{j}|^2}{c}-\lambda_0\frac{|\varrho_0|^2}{c}\\
&\leq\frac{1}{C_{*}}\zeta^0+\frac{1}{C_{*}}\delta(r)c-\lambda_0\frac{|\varrho_0|^2}{c}\\
&\leq\lambda_0\left(\frac{\zeta^0}{C_{*}\lambda_0}-\frac{|\varrho_0|^2}{c}\right)+\frac{1}{C_{*}}\delta(r)c.
\end{align*}
It follows that $$\sup\limits_{u\in \mathcal{M}_{c}^{r}}\frac{\|u-\varrho_0\psi_0\|^2_{\dot{\Sigma}}}{c}
\leq\left(\frac{\lambda_0}{\inf\limits_{j\geq1}(\lambda_{j}
-\lambda_0)}+1\right)\frac{2}{C_{*}}\delta(r)c,
$$
that is, \eqref{11} holds. This completes the proof.
\qed

In the following, we prove Theorem \ref{thm1'}.

\noindent\textbf{Proof of Theorem \ref{thm1'}.}
We first claim that $u_{c}\to0$ in $\Sigma$ as $c\to0^{+}$.
Since $(u_{c}, \mu_{c})\in\mathcal{M}_{c}^{r}\times\mathbb{R}^{-}$ is a couple solution to \eqref{EQ}, we have
\begin{align}\label{Q2}
\frac{1}{2}\int_{\mathbb{R}^3}|\nabla u_{c}|^2dx+\int_{\mathbb{R}^3}V(x)|u_{c}|^2dx+N(u_{c})
-\Omega\int_{\mathbb{R}^3}\overline{u}_{c}L_{z}u_{c}dx+\mu_{c}\int_{\mathbb{R}^3}|u_{c}|^2dx=0.
\end{align}
By Lemma \ref{lem5} and $\mu_{c}<0$ with $0<c<c_0$, the energy functional $E_{\Omega}(u)$ can be rewritten as
\begin{align*}
E_{\Omega}(u_{c})&=\frac{1}{4}\int_{\mathbb{R}^3}|\nabla u_{c}|^2dx+\frac{1}{2}\int_{\mathbb{R}^3}V(x)|u_{c}|^2dx-\frac{1}{2}\Omega\int_{\mathbb{R}^3}\overline{u}_{c}L_{z}u_{c}dx-\frac{\mu_{c}}{2}\int_{\mathbb{R}^3}|u_{c}|^2dx\\
&\geq\frac{C_{*}}{4}\int_{\mathbb{R}^3}|\nabla u_{c}|^2+|x|^2|u|^2dx.
\end{align*}
It follows from \eqref{SJ} that
$$
E_{\Omega}(u_{c})=m(c,r)\leq\zeta^0c+\sqrt{2}\Lambda \mathcal{C}_4^4(\zeta^0)^{\frac{3}{2}}c^2,
$$
where $\zeta^0=\frac{1}{2}\sum\limits_{j=1}\limits^3\gamma_{j}$.
It yields that
$$
\frac{C_{*}}{4}\int_{\mathbb{R}^3}|\nabla u_{c}|^2+|x|^2|u_{c}|^2dx\leq \zeta^0c+\sqrt{2}\Lambda \mathcal{C}_4^4(\zeta^0)^{\frac{3}{2}}c^2\to0
$$
as $c\to0^{+}$. Namely, $u_{c}\to0$ in $\Sigma$ as $c\to0^{+}$.

According to \eqref{mu-c-bound}, we know that
$$-\zeta^0-\delta(r)c\leq\mu_{c}\leq3\left[-\frac{C_{*}}{2}+\Lambda \mathcal C_4^4r^{\frac{1}{2}}c^{\frac{1}{2}}\right]$$
for $0<c<c_0$, where $\zeta^0=\frac{1}{2}\sum\limits_{j=1}\limits^{3}\gamma_{j}$ and $\delta(r)=\sqrt{2}\Lambda \mathcal{C}_4^4(\zeta^0)^{\frac{3}{2}}+\frac{\Lambda}{2}\mathcal C_4^4r^{\frac{3}{2}}$.
Then there exists a constant $\mu_0\in[-\zeta^0, -\frac{3C_{*}}{2}]$ such that $\mu_{c}\to \mu_0$ as $c\to0^{+}$.
By Lemma \ref{lem2}, \eqref{GN} and the fact that $u_{c}\in B(cr)$, we have
\begin{align}\label{K}
0<\frac{|N(u_{c})|}{|u_{c}|_2^2}\leq\frac{\Lambda |u_{c}|_4^4}{|u_{c}|_2^2}\leq\frac{\Lambda \mathcal C_4^4|\nabla u_{c}|_2^3}{|u_{c}|_2}\leq\Lambda \mathcal C_4^4r^{\frac{3}{2}}c\to0
\end{align}
as $c\to0^{+}$.
Note that $u_c$ satisfies the Pohozaev identity
\eqref{D5}, i.e., $Q(u_c)=0$.
Set $w_{c}:=\frac{u_{c}}{|u_{c}|_2}$, then we can deduce that
\begin{align*}
0=\frac{Q(u_{c})}{|u_{c}|_2^2}=\frac{1}{2}\int_{\mathbb{R}^3}|\nabla w_{c}|^2dx-\int_{\mathbb{R}^3}V(x)|w_{c}|^2dx+\frac{3}{4}\frac{N(u_{c})}{|u_{c}|_2^2},
\end{align*}
which yields that $\lim\limits_{c\to0^{+}}\int_{\mathbb{R}^3}|\nabla w_{c}|^2dx=\lim\limits_{c\to0^{+}}\int_{\mathbb{R}^3}2V(x)|w_{c}|^2dx$.
From \eqref{Q2} and \eqref{K},
we infer that
\begin{align*}
\lim\limits_{c\to0^{+}}\mu_{c}&=
\lim\limits_{c\to0^{+}}\left[-\frac{1}{2}\int_{\mathbb{R}^3}|\nabla w_{c}|^2+2V(x)|w_{c}|^2dx
+\int_{\mathbb{R}^3}\Omega\overline{w}_{c}L_{z}w_{c}dx-\frac{N(u_{c})}{|u_{c}|_2^2}\right]\\
&=\lim\limits_{c\to0^{+}}\left[-\frac{1}{2}\int_{\mathbb{R}^3}|\nabla w_{c}|^2+2V(x)|w_{c}|^2dx
+\int_{\mathbb{R}^3}\Omega\overline{w}_{c}L_{z}w_{c}dx\right].
\end{align*}
Similarly, by \eqref{K}, one has
\begin{align*}
\lim\limits_{c\to0^{+}}\frac{m(c,r)}{c}&=\lim\limits_{c\to0^{+}}\left[\frac{1}{2}\int_{\mathbb{R}^3}|\nabla w_{c}|^2+2V(x)|w_{c}|^2dx
-\int_{\mathbb{R}^3}\Omega\overline{w}_{c}L_{z}w_{c}dx-\frac{1}{2}\frac{N(u_{c})}{|u_{c}|_2^2}\right]\\
&=\lim\limits_{c\to0^{+}}\left[\frac{1}{2}\int_{\mathbb{R}^3}|\nabla w_{c}|^2+2V(x)|w_{c}|^2dx
-\int_{\mathbb{R}^3}\Omega\overline{w}_{c}L_{z}w_{c}dx\right].
\end{align*}
Therefore, we conclude that
$$
\lim\limits_{c\to0^{+}}\frac{m(c,r)}{c}=
\lim\limits_{c\to0^{+}}-\mu_{c}=-\mu_0
$$
 and
\begin{align*}
\lim\limits_{c\to0^{+}}\frac{\int_{\mathbb{R}^3}|\nabla u_{c}|^2dx
-\int_{\mathbb{R}^3}\Omega\overline{u}_{c}L_{z}u_{c}dx}{c}=\lim\limits_{c\to0^{+}}\frac{\int_{\mathbb{R}^3}2V(x)|u_{c}|^2dx
-\int_{\mathbb{R}^3}\Omega\overline{u}_{c}L_{z}u_{c}dx}{c}=-\mu_0.
\end{align*}

Finally, it remains to prove the local minimizer $u_c$ is a ground state
if $0<\Omega<\frac{\sqrt{5}}{3}\min\{\gamma_1, \gamma_2\}$ and $c>0$ sufficiently small. Suppose by contradiction that there exists $v\in S(c)$ such that
\begin{align*}
(E_{\Omega}|_{S(c)})'(v)=0\,\,\,\,{\rm{and}}\,\,\,\,E_{\Omega}(v)<m(c,r).
\end{align*}
Then we can infer that $Q(v)=0$, where
\begin{align*}
Q(v)=\frac{1}{2}\int_{\mathbb{R}^3}|\nabla v|^2dx-\int_{\mathbb{R}^3}V(x)|v|^2dx+\frac{3}{4}\int_{\mathbb{R}^3}\lambda_1|v|^4+\lambda_2(K\ast|v|^2)|v|^2dx.
\end{align*}
Similar to the proof of \eqref{bounded-PS-sequence}, we conclude that
\begin{align}\label{P1}
E_{\Omega}(v)\geq\widetilde{C}\int_{\mathbb{R}^3}|\nabla v|^2+|x|^2|v|^2dx.
\end{align}
Combining this with \eqref{SJ}, we obtain that
\begin{align*}
\widetilde{C}\int_{\mathbb{R}^3}|\nabla v|^2+|x|^2|v|^2dx\leq E_{\Omega}(v)<m(c,r)\leq \xi^0c+\sqrt{2}\Lambda \mathcal{C}_4^4(\xi^0)^{\frac{3}{2}}c^2\to0,
\end{align*}
as $c\to0$. This yields that $v\in B(r)$ for $c$ small enough. By the definition of $m(c,r)$, we have $E_{\Omega}(v)\geq m(c,r)$,
which contradicts $E_{\Omega}(v)<m(c,r)$.
This completes the proof.\qed

\section{The critical rotational speed }\label{critical}
In this section, we mainly consider the case of $0<\Omega=\Omega^{*}$ (the critical rotational speed) and
investigate the existence of normalized solutions for \eqref{EQ}. Namely, the proof of Theorem \ref{thm3} is given.

Since $(\lambda_1, \lambda_2)\in \mathbb{R}^2$ satisfies \eqref{UR'}, the dipolar term $N(u)<0$, which is $L^2$-supercritical.
A standard scaling argument implies that the energy functional $E_{\Omega^{*}}(u)$ is no longer bounded from below on $S(c)$.
Motivated by \cite{Bellazzini2017}, we study the following local minimization problem: for any given $\rho>0$,
\begin{align}\label{4.4}
m_{\Omega^{*}}(c,\rho):=\inf\{E_{\Omega^{*}}(u): u\in S(c)\cap B_{\Omega^{*}}(\rho)\},\,\,\,\,c>0,
\end{align}
where $S(c)=\{u\in \Sigma_{\Omega^*}:|u|_2^2=c\}$ and
\begin{align}\label{B-rho-definition}
B_{\Omega^{*}}(\rho):=\left\{u\in \Sigma_{\Omega^{*}}: \int_{\mathbb{R}^3}|(\nabla-iA)u|^2+2V_{\Omega^*}(x)|u|^2dx\leq \rho\right\}.
\end{align}
As we will see, the existence of minimizers of $m_{\Omega^{*}}(c,\rho)$ can be established by applying the concentration-compactness principle.
We first introduce the following lemma, which plays a important role in excluding the vanishing scenario.

\begin{lem}\label{lem15}
Let $\{u_{n}\}$ be a sequence of $E_{\Omega^{*}}$ in $\Sigma_{\Omega^{*}}$. Suppose that
\begin{align}\label{4.5}
\sup\limits_{n\geq1}\|u_{n}\|_{\Sigma_{\Omega^{*}}}<\infty,
\end{align}
and there is $\varepsilon_0>0$ such that
\begin{align}\label{4.6}
\inf\limits_{n\geq1}\int_{\mathbb{R}^3}|u_{n}|^4dx\geq\varepsilon_0>0.
\end{align}
Then there exist $u\in \Sigma_{\Omega^{*}}\backslash\{0\}$ and $\{z_{n}\}\subseteq\mathbb{R}^3$ with $z_{n}=(0, z_{n}^2, 0)$ such that, up to a subsequence,
\begin{align*}
w_{n}(x):=e^{iA(z_{n})\cdot x}u_{n}(x+z_{n})\rightharpoonup u\,\,\,\,{\rm{in}}\,\,\,\,\Sigma_{\Omega^{*}}.
\end{align*}
\end{lem}
\begin{prf}
By the H\"{o}lder's inequality and Sobolev embedding theorem, we obtain that
\begin{align*}
\int_{T_{k}}|u_{n}|^4dx&\leq \left(\int_{T_{k}}|u_{n}|^2dx\right)^{\frac{1}{2}}\left(\int_{T_{k}}|u_{n}|^6dx\right)^{\frac{1}{2}}\\
&\leq C|u_{n}|_{L^2(T_{k})}\left(\int_{T_{k}}|\nabla|u_{n}||^2+|u_{n}|^2dx\right)^{\frac{3}{2}},
\end{align*}
where
\begin{align*}
T_{k}:=\mathbb{R}\times(k,k+1)\times\mathbb{R},\,\,\,k\in\mathbb{Z}.
\end{align*}
Taking a sum over $k\in \mathbb{Z}$, it follows from the diamagnetic inequality \eqref{DI} that
\begin{align}\label{4.7}
\int_{\mathbb{R}^3}|u_{n}|^4dx&\leq
C\left(\sup\limits_{k\in\mathbb{Z}}|u_{n}|_{L^2(T_{k})}
\right)\left(\int_{\mathbb{R}^3}|\nabla|u_{n}||^2+|u_{n}|
^2dx\right)^{\frac{3}{2}}\notag\\
&\leq
C\left(\sup\limits_{k\in\mathbb{Z}}|u_{n}|_{L^2(T_{k})}
\right)\left(\int_{\mathbb{R}^3}|(\nabla-iA)u_{n}|^2+|u_{n}|^2dx
\right)^{\frac{3}{2}}.
\end{align}
Taking account of \eqref{4.5}-\eqref{4.6}, we derive from \eqref{4.7} that there is $\{k_{n}\}\subseteq\mathbb{Z}$ such that
\begin{align}\label{4.8}
\inf\limits_{n\geq1}\int_{T_{k_{n}}}|u_{n}|^2dx\geq C>0
\end{align}
for some constant $C>0$. Next, we can prove that up to a subsequence if necessary,
there exist $u\in \Sigma_{\Omega^{*}}\backslash\{0\}$ and $\xi_{n}:=(0, k_{n}, 0)$ with $k_{n}\in\mathbb{Z}$ such that
\begin{align}\label{A1}
e^{iA(\xi_{n})\cdot x}u_{n}(x+\xi_{n})\rightharpoonup u\,\,\,\,{\rm{in}}\,\,\,\,\Sigma_{\Omega^{*}},
\end{align}
where $A(\xi_{n})=\Omega^*(-k_{n},0,0)$.
In fact, by \eqref{4.5} and \eqref{4.8}, we deduce that
\begin{align}\label{4.9}
\int_{T_0}|e^{iA(\xi_{n})\cdot x}u_{n}(x+\xi_{n})|^2dx=\int_{T_{k_{n}}}|u_{n}|^2dx\geq C>0,
\end{align}
and
$$\small
\begin{aligned}
&\sup\limits_{n\geq1}\int_{T_0}|(\nabla- iA)u_{n}(x+\xi_{n})e^{iA(\xi_{n})\cdot x}|^2
+2V_{\Omega^{*}}(x)|u_{n}(x+\xi_{n})e^{iA(\xi_{n})\cdot x}|^2+|u_{n}(x+\xi_{n})e^{iA(\xi_{n})\cdot x}|^2dx\\
\,\,\,&=\sup\limits_{n\geq1}\int_{T_{k_{n}}}|(\nabla-iA)u_{n}|^2
+2V_{\Omega^{*}}(x)|u_{n}(x)|^2+|u_{n}(x)|^2dx<+\infty,
\end{aligned}
$$
where $T_{k_{n}}=\mathbb{R}\times(k_{n}, k_{n}+1)\times\mathbb{R}$, $k_{n}\in\mathbb{Z}$ and $T_0=\mathbb{R}\times(0,1)\times\mathbb{R}$.
Hence up to a subsequence if necessary, there exists a $u\in \Sigma_{\Omega^{*}}$ such that
\begin{align*}
e^{iA(\xi_{n})\cdot x}u_{n}(x+\xi_{n})\rightharpoonup u\,\,\,{\rm{in}}\,\,\,\Sigma_{\Omega^{*}}\,\,\,{\rm{as}}\,\,\,n\to\infty.
\end{align*}
By the Sobolev compact embedding theorem, we have that
\begin{align}\label{4.10}
e^{iA(\xi_{n})\cdot x}u_{n}(x+\xi_{n})\to u\,\,\,{\rm{in}}\,\,\,L_{loc}^2(T_0)\,\,\,{\rm{as}}\,\,\,n\to\infty,
\end{align}
where $T_0=\mathbb{R}\times(0,1)\times\mathbb{R}$. Note that for any $\varepsilon>0$, there exists a constant $R>0$ such that
\begin{align}\label{4.11}
\frac{1}{|x|^2}\leq\varepsilon\,\,\,{\rm{for}}\,\,\,|x|\geq R.
\end{align}
From \eqref{4.10}, we obtain that for the fixed $R>0$,
\begin{align*}
e^{iA(\xi_{n})\cdot x}u_{n}(x+\xi_{n})\to u\,\,\,{\rm{in}}\,\,\, L^2(T_0\cap B_{R}(0))\,\,\,{\rm{as}}\,\,\,n\to\infty.
\end{align*}
It yields that there is a positive integer $N$ such that
\begin{align*}
\int_{T_0\cap B_{R}(0)}|e^{iA(\xi_{n})\cdot x}u_{n}(x+\xi_{n})-u|^2dx\leq\varepsilon\,\,\,{\rm{for}}\,\,\,n\geq N.
\end{align*}
Combining \eqref{4.5} and \eqref{4.11}, we conclude that for $n\geq N$,
\begin{align*}
\int_{T_0}|e^{iA(\xi_{n})\cdot x}u_{n}(x+\xi_{n})-u|^2dx
=&\int_{T_0\cap B_{R}(0)}|e^{iA(\xi_{n})\cdot x}u_{n}(x+\xi_{n})-u|^2dx\\
&+\int_{T_0\setminus (T_0\cap B_{R}(0))}|e^{iA(\xi_{n})\cdot x}u_{n}(x+\xi_{n})-u|^2dx\\
\leq&\varepsilon+\varepsilon\int_{T_0\setminus (T_0\cap B_{R}(0))}|x|^2|e^{iA(\xi_{n})\cdot x}u_{n}(x+\xi_{n})-u|^2dx\\
\leq&\varepsilon+\varepsilon\int_{T_0\setminus (T_0\cap B_{R}(0))}(x_1^2+1+x_3^3)|e^{iA(\xi_{n})\cdot x}u_{n}(x+\xi_{n})-u|^2dx\\
\leq& C\varepsilon,
\end{align*}
which implies that $e^{iA(\xi_{n})\cdot x}u_{n}(x+\xi_{n})\to u$ in $L^2(T_0)$. In view of \eqref{4.9}, one has $u\neq0$, which implies that the claim \eqref{A1} holds.
This completes the proof.
\end{prf}

In the following, we introduce the $L^2$-bound of the magnetic Sobolev norm, which comes from \cite[Proposition 2.2 and Remark 2.5]{Esteban1989}.

\begin{prop}\label{prop1}
Let $A=(A_1,..., A_{N})\in W^{1,\infty}_{loc}(\mathbb{R}^{N}, \mathbb{R}^{N})$ and $j,k\in\{1,...,N\}$. Then for any $u\in C_0^{\infty}(\mathbb{R}^{N})$, we have
\begin{align*}
\left|\int_{\mathbb{R}^{N}}(\partial_{j}A_{k}-\partial_{k}A_{j})u\overline{u}dx\right|\leq |(\partial_{j}-iA_{j})u|_2^2+|(\partial_{k}-iA_{k})u|_2^2.
\end{align*}
In particular, if $N=3$ and $A(x)=\Omega^*(-x_2, x_1,0)$, then for any $u\in C_0^{\infty}(\mathbb{R}^3)$, we have
\begin{align}\label{MS}
2\Omega^*|u|_2^2\leq |(\partial_{1}-iA_{1})u|_2^2+|(\partial_{2}-iA_{2})u|_2^2.
\end{align}
Moreover, if $N=2$, then \eqref{MS} is achieved by $u(x)=\sqrt{\frac{\Omega^*}{\pi}}e^{-\frac{\Omega^*}{2}|x|^2}$.
\end{prop}

In order to apply the Lemma \ref{lem15} to guarantee the non-vanishing of the minimizing sequences up to translation, we need to prove the following result.
\begin{lem}\label{lem16}
Let $\{u_{n}\}$ be a minimizing sequence for $m_{\Omega^{*}}(c,\rho)$. Then for any given $\rho>0$,
there exist constants $\delta>0$ and $c_0(\rho)>0$ depending on $\rho$ such that
\begin{align*}
\liminf\limits_{n\to\infty}|N(u_{n})|>\delta>0\,\,\,\,\,{\rm{for}}\,\,\,\,{\rm{all}}\,\,\,0<c<c_0(\rho).
\end{align*}
\end{lem}
\begin{prf}
We assume by contradiction that the minimizing sequence $\{u_{n}\}$, up to a subsequence, satisfies that $\lim\limits_{n\to\infty}|N(u_{n})|=0$.
Denote $x=(x_1,x_2,x_3)\in\mathbb{R}^3$ and $A(x)=(A_1,A_2,A_3)$. By the definition of $A(x)$,
we know that $A_1=-\Omega^{*}x_2$, $A_2=\Omega^{*}x_1$ and $A_3=0$.
Set
\begin{align}\label{S1}
\zeta_{\Omega^{*}}^{0}:=\inf\left\{\frac{1}{2}\int_{\mathbb{R}^3}|(\nabla-iA) u|^2dx+\int_{\mathbb{R}^3}V_{\Omega^{*}}(x)|u|^2dx: u\in \Sigma_{\Omega^{*}}, |u|_2^2=1\right\},
\end{align}
where $V_{\Omega^{*}}(x)=\frac{1}{2}(\gamma_1^2-\gamma_2^2)x_1^2+\frac{1}{2}\gamma_3^2x_3^2$ with $\Omega^{*}=\gamma_2$.
Inspired by \cite[Lemma 2.1]{Bellazzini2017}, we first claim that
\begin{align}\label{S2}
\zeta_{\Omega^{*}}^{0}\geq\gamma_2+\frac{\gamma_3}{2}.
\end{align}
To this end,
we introduce a Hilbert basis $\Phi_{k}(x_3)$ of $L^2(\mathbb{R})$ and $\lambda_{k}$ for $k\geq0$ such that
\begin{align*}
-\frac{1}{2}\Delta_{x_3}\Phi_{k}+\frac{1}{2}\gamma_3^2x_3^2\Phi_{k}&=\lambda_{k}\Phi_{k},\,\,\,\,\int_{\mathbb{R}}|\Phi_{k}|^2dx_3=1,\\
\lambda_{k}&\leq\lambda_{k+1},\,\,\, k=0,1,2,\cdots.
\end{align*}
It follows from \eqref{T} that $\lambda_0=\frac{\gamma_3}{2}$.
By the Fourier decomposition, we have
\begin{align}\label{33}
v(x)=\sum\limits_{k\geq0}v_{k}(x_1,x_2)\Phi_{k}(x_3).
\end{align}
Let $|v|_2^2=1$, then
\begin{align}\label{32}
1=\int_{\mathbb{R}^3}|v(x)|^2dx&=\sum\limits_{k\geq0}\left(\int_{\mathbb{R}^2}|v_{k}|^2dx_1dx_2\right)\left(\int_{\mathbb{R}}|\Phi_{k}|^2dx_3\right)\notag\\
&=\sum\limits_{k\geq0}\int_{\mathbb{R}^2}|v_{k}|^2dx_1dx_2.
\end{align}
From \eqref{33}, \eqref{32} and Proposition \ref{prop1}, we deduce that
\begin{align*}
\frac{1}{2}\int_{\mathbb{R}^3}&|(\nabla-iA) v|^2dx+\int_{\mathbb{R}^3}V_{\Omega^{*}}(x)|v|^2dx\\
=&\frac{1}{2}\sum\limits_{j=1}^{2}
\int_{\mathbb{R}^3}\left|\left(\frac{\partial}{\partial x_{j}}-iA_{j}\right)v\right|^2dx
+\frac{1}{2}\int_{\mathbb{R}^3}\left|\frac{\partial v}{\partial x_{3}}\right|^2dx+\int_{\mathbb{R}^3}V_{\Omega^{*}}|v|^2dx\\
=&\frac{1}{2}\sum\limits_{k\geq0}\sum\limits_{j=1}^{2}\int_{\mathbb{R}^2}\left|\left(\frac{\partial}{\partial x_{j}}-iA_{j}\right)v_{k}\right|^2dx_1dx_2
+\frac{1}{2}\sum\limits_{k\geq0}\int_{\mathbb{R}}\left|\frac{\partial \Phi_{k} }{\partial x_{3}}\right|^2+\gamma_3^2x_3^2|\Phi_{k}|^2dx_3\\
&+\frac{1}{2}\int_{\mathbb{R}^3}(\gamma_1^2-\gamma_2^2)x_1^2|v|^2dx\\
\geq& \Omega^{*}+\frac{1}{2}\sum\limits_{k\geq0}\int_{\mathbb{R}}\lambda_{k}|\Phi_{k}|^2dx_3\\
\geq&\Omega^{*}+\frac{1}{2}\lambda_0\sum\limits_{k\geq0}\int_{\mathbb{R}}|\Phi_{k}|^2dx_3=\gamma_2+\frac{\gamma_3}{2},
\end{align*}
which yields that \eqref{S2} holds.
It follows  that
\begin{align}\label{4.12}
m_{\Omega^{*}}(c,\rho)&=E_{\Omega^{*}}(u_{n})+o_{n}(1)\notag\\
&=\frac{1}{2}\int_{\mathbb{R}^3}|(\nabla-iA)u_{n}|^2dx+\int_{\mathbb{R}^3} V_{\Omega^{*}}(x)|u_{n}|^2dx+\frac{1}{2}N(u_{n})+o_{n}(1)\notag\\
&=\frac{1}{2}\int_{\mathbb{R}^3}|(\nabla-iA)u_{n}|^2dx+\int_{\mathbb{R}^3} V_{\Omega^{*}}(x)|u_{n}|^2dx+o_{n}(1)\notag\\
&\geq\zeta_{\Omega^{*}}^0|u_{n}|_2^2++o_{n}(1)\geq(\gamma_2+\frac{\gamma_3}{2})c.
\end{align}
Next, in order to establish a contradiction,
we take $u_0(x_1,x_2):=\sqrt{\frac{\Omega^{*}}{\pi}}e^{-\frac{\Omega^{*}}{2}(x_1^2+x_2^2)}$.
 Clearly, $|u_0(x_1,x_2)|_2^2=1$ and
\begin{align}\label{4.13}
2\Omega^*= |(\partial_{1}-iA_{1})u_0(x_1,x_2)|_2^2+|(\partial_{2}-iA_{2})u_0(x_1,x_2)|_2^2
\end{align}
due to Proposition \ref{prop1}. Set $\varphi_{c}(x_3):=\sqrt{c}\Phi_0(x_3)$, where $\Phi_0(x_3)\in \Sigma(\mathbb{R})$ is as in \eqref{E},
that is, $\Phi_0(x_3)=\left(\frac{\gamma_3}{\pi}\right)^{\frac{1}{4}}e^{-\frac{\gamma_3}{2}x_3^2}$. Then the direct calculation implies that
\begin{align}\label{4.14}
|\varphi_{c}(x_3)|_2^2=c,\,\,\,\,\frac{1}{2}\int_{\mathbb{R}}|\nabla\varphi_{c}(x_3)|^2+\gamma_3^2x_3^2|\varphi_{c}(x_3)|^2dx_3=\frac{c}{2}\gamma_3.
\end{align}
Now we define $v_{c}(x)=u_0(x_1,x_2)\varphi_{c}(x_3)$.
Then it is obvious that there exists $c_0(\rho)>0$ such that
$v_{c}(x)\in B_{\Omega^{*}}(\rho)$ provided that $c<c_0(\rho)$, where $c_0(\rho)$ is a constant depending on $\rho$.
So $v_{c}(x)\in S(c)\cap B_{\Omega^{*}}(\rho)$. By \eqref{4.13} and \eqref{4.14}, we further deduce that
\begin{align*}
m_{\Omega^{*}}(c,\rho)\leq&E_{\Omega^{*}}(v_{c}(x))\\
=&\frac{1}{2}\int_{\mathbb{R}^3}|(\nabla-iA)v_{c}|^2dx
+\int_{\mathbb{R}^3}V_{\Omega^*}(x)|v_{c}|^2dx
+\frac{1}{2}N(v_{c})\\
=&\frac{1}{2}\int_{\mathbb{R}^3}\left|\left(\frac{\partial}{\partial{x_1}}-iA_1\right)v_{c}\right|^2
+\left|\left(\frac{\partial}{\partial{x_2}}-iA_2\right)v_{c}
\right|^2+\left|\frac{\partial v_{c}}{\partial x_3}\right|^2dx\\
& +\int_{\mathbb{R}^3}V_{\Omega^*}(x)|v_{c}|^2dx+\frac{1}{2}N(v_{c})\\
=&\frac{1}{2}\int_{\mathbb{R}^2}\left|\left(\frac{\partial}{\partial{x_1}}-iA_1\right)u_0(x_1,x_2)\right|^2dx_1dx_2\int_{\mathbb{R}}|\varphi_{c}(x_3)|^2dx_3\\&+
\frac{1}{2}\int_{\mathbb{R}^2}\left|\left(\frac{\partial}{\partial{x_2}}-iA_2\right)u_0(x_1,x_2)\right|^2dx_1dx_2\int_{\mathbb{R}}\left|\varphi_{c}(x_3)\right|^2dx_3\\&+
\frac{1}{2}\int_{\mathbb{R}}\left|\frac{\partial \varphi_{c}(x_3)}{\partial x_3}\right|^2dx_3\int_{\mathbb{R}^2}|u_0(x_1,x_2)|^2dx_1dx_2\\
&+\frac{1}{2}\int_{\mathbb{R}}\gamma_3^2x_3^2|\varphi_{c}(x_3)|^2
dx_3\int_{\mathbb{R}^2}|u_0(x_1,x_2)|^2dx_1dx_2\\
& +\int_{\mathbb{R}^3}(V_{\Omega^*}(x)-\frac{1}{2}\gamma_3^2x_3^2)|v_{c}|^2dx+\frac{1}{2}N(v_{c})\\
=&\gamma_2c+\frac{c}{2}\gamma_3+\int_{\mathbb{R}^3}\frac{1}{2}(\gamma_1^2-\gamma_2^2)x_1^2|v_{c}|^2dx+\frac{1}{2}N(v_{c})\\
=&\gamma_2c+\frac{c}{2}\gamma_3+\int_{\mathbb{R}^2}\frac{1}{2}
(\gamma_1^2-\gamma_2^2)x_1^2|u_0(x_1,x_2)|^2dx\int_{\mathbb{R}}\varphi_{c}(x_3)|^2dx+\frac{1}{2}N(v_{c})\\
=&\gamma_2c+\frac{c}{2}\gamma_3+\frac{(\gamma_1^2-\gamma_2^2)c}{4\gamma_2}+\frac{1}{2}N(v_{c}).
\end{align*}
By applying the Fourier transform for $N(v_c)$ with $(\lambda_1,\lambda_2)$ satisfying \eqref{UR'}, we obtain that
\begin{align*}
N(v_{c})=&\int_{\mathbb{R}^3}\lambda_1|v_{c}|^4+\lambda_2(K*|v_{c}|^2)|v_{c}|^2dx\\
=&\frac{1}{(2\pi)^3}\int_{\mathbb{R}^3}[\lambda_1+\lambda_2\widehat{K}(\xi)]|\widehat{v_{c}^2}|^2d\xi\\
=&\frac{1}{(2\pi)^3}\int_{\mathbb{R}^3}\left[\lambda_1+\frac{4\pi}{3}\lambda_2\left(\frac{3\xi_3^2}{|\xi|^2}-1\right)\right]|\widehat{v_{c}^2}|^2d\xi\\
=&\frac{1}{(2\pi)^3}\int_{\mathbb{R}^3}\left[\lambda_1+\frac{4\pi}{3}\lambda_2\left(\frac{3\xi_3^2}{|\xi|^2}-1\right)\right]
\left|\int_{\mathbb{R}^3}e^{-i\cdot\xi}v_{c}^2(x)dx\right|^2d\xi\\
=&\frac{1}{(2\pi)^3}\int_{\mathbb{R}^3}\left[\lambda_1+\frac{4\pi}{3}\lambda_2\left(\frac{3\xi_3^2}{|\xi|^2}-1\right)\right]
\left|\frac{\gamma_2\gamma_3^{\frac{1}{2}}c}{\pi^{\frac{3}{2}}}\int_{\mathbb{R}^3}e^{-i\cdot\xi}\cdot e^{-\gamma_2(x_1^2+x_2^2)-\gamma_3x_3^2}dx\right|^2d\xi\\
=&-\frac{\gamma_2^2\gamma_3c^2}{8\pi^6}\int_{\mathbb{R}^3}\left|\lambda_1+\frac{4\pi}{3}\lambda_2\left(\frac{3\xi_3^2}{|\xi|^2}-1\right)\right|
\left|\int_{\mathbb{R}^3}e^{-i\cdot\xi}\cdot e^{-\gamma_2(x_1^2+x_2^2)-\gamma_3x_3^2}dx\right|^2d\xi\\
\triangleq&-M(\gamma_2, \gamma_3),
\end{align*}
where $\widehat{u}(\xi):=\mathcal{F}u(\xi)
=\int_{\mathbb{R}^3}e^{-ix\cdot\xi}u(x)dx$
and
$M(\gamma_2, \gamma_3)$ is a constant depending on $\gamma_2$, $\gamma_3$.
Choosing $\gamma_2\leq\gamma_1
<\sqrt{\frac{2M(\gamma_2,\gamma_3)\gamma_2}{c}+\gamma_2^2}$, we can deduce that $\frac{(\gamma_1^2-\gamma_2^2)c}{4\gamma_2}+\frac{1}{2}N(v_{c})<0$, and then
$$
m_{\Omega^{*}}(c,\rho)<(\gamma_2+\frac{\gamma_3}{2})c,
$$
which contradicts \eqref{4.12}. This completes the proof.\\
\end{prf}

\noindent\textbf{Proof of Theorem \ref{thm3}}.
$(i)$ Since $N(u)<0$, similar to the proof of \cite[Lemma 3.1]{Bellazzini2017}, we conclude that for any given $\rho>0$ there is $0<c_1(\rho)<1$ such that
\begin{align}\label{V1}
S(c)\cap B_{\Omega^{*}}(\rho)\neq\emptyset,
\end{align}
\begin{align}\label{V2}
\inf\limits_{u\in S(c)\cap B_{\Omega^{*}}(\rho\frac{c}{2})}E_{\Omega^{*}}(u)
<\inf\limits_{u\in S(c)\cap \left(B_{\Omega^{*}}(\rho)\setminus B_{\Omega^{*}}(\rho c)\right)}E_{\Omega^{*}}(u),
\end{align}
for any $0<c<c_1(\rho)$. Thus, $m_{\Omega^{*}}(\rho, c)$ is well-defined for $0<c<c_1(\rho)$.
Let $\{u_{n}\}\subseteq S(c)\cap B_{\Omega^{*}}(\rho)$ be a minimizing sequence for $m_{\Omega^{*}}(\rho, c)$
satisfying $\lim\limits_{n\to\infty}E_{\Omega^{*}}(u_{n})=m_{\Omega^{*}}(\rho, c)$. It is clear that $\{u_{n}\}$ is bounded in $\Sigma_{\Omega^{*}}$.
By Lemma \ref{lem16}, we take a subsequence still denote by $\{u_{n}\}$ such that
\begin{align*}
\liminf\limits_{n\to\infty}|N(u_{n})|>\delta>0
\end{align*}
for some constant $\delta>0$. Combining this with Lemma \ref{lem15},
there exist $\widetilde{u}_{c}\in \Sigma_{\Omega^{*}}\backslash\{0\}$ and $\{z_{n}\}\subseteq\mathbb{R}^3$ with $z_{n}=(0, z_{n}^2, 0)$ such that, up to a subsequence,
\begin{align*}
w_{n}(x):=e^{iA(z_{n})\cdot x}u_{n}(x+z_{n})\rightharpoonup\widetilde{u}_{c} \,\,\,\,{\rm{in}}\,\,\,\,\Sigma_{\Omega^{*}}.
\end{align*}
It is obvious that
\begin{align*}
\int_{\mathbb{R}^3}|w_{n}(x)|^2dx&=\int_{\mathbb{R}^3}|u_{n}|^2dx
\ \mbox{and}\
N(w_{n})=N(u_{n}).
\end{align*}
Moreover, we know that
\begin{align*}
\int_{\mathbb{R}^3}|(\nabla-iA)w_{n}|^2+V_{\Omega^{*}}(x)|w_{n}|^2dx&=\int_{\mathbb{R}^3}|(\nabla-iA)u_{n}|^2+V_{\Omega^{*}}(x)|u_{n}|^2dx,
\end{align*}
see \cite[Eq.(3.16)]{Esteban1989}.
Thus, $E_{\Omega^{*}}(w_{n})=E_{\Omega^{*}}(u_{n})$. That is to say, $\{w_{n}\}$ is also a minimizing sequence of $m_{\Omega^{*}}(c,\rho)$. By the weakly lower
semi-continuous of norm, we get that
\begin{align*}
0<\int_{\mathbb{R}^3}|\widetilde{u}_{c}|^2dx\leq\liminf\limits_{n\to\infty}\int_{\mathbb{R}^3}|w_{n}|^2dx=\liminf\limits_{n\to\infty}\int_{\mathbb{R}^3}|u_{n}|^2dx=c
\end{align*}
and
\begin{align*}
\int_{\mathbb{R}^3}|(\nabla-iA)\widetilde{u}_{c}|^2+2V_{\Omega^{*}}(x)|\widetilde{u}_{c}|^2dx&\leq\int_{\mathbb{R}^3}|(\nabla-iA)w_{n}|^2+2V_{\Omega^{*}}(x)|w_{n}|^2dx\\
&=\int_{\mathbb{R}^3}|(\nabla-iA)u_{n}|^2+V_{\Omega^{*}}(x)|u_{n}|^2dx\\
&\leq\rho,
\end{align*}
which imply that $\widetilde{u}_{c}\in B_{\Omega^{*}}(\rho)$. To claim that $\widetilde{u}_{c}$ is a minimizer for $m_{\Omega^{*}}(c,\rho)$. We first show that $|\widetilde{u}_{c}|_2^2=c$.
Assume that $0<|\widetilde{u}_{c}|_2^2=r<c$. By \cite[Lemma 2.2]{Wu2024}, one has
\begin{align*}
N(w_{n})=N(\widetilde{u}_{c})+N(w_{n}-\widetilde{u}_{c})+o_{n}(1).
\end{align*}
By the weak convergence of $w_{n}$ in $\Sigma_{\Omega^{*}}$, we deduce that
\begin{align*}
\int_{\mathbb{R}^3}|(\nabla-iA)w_{n}|^2dx&=\int_{\mathbb{R}^3}|(\nabla-iA)\widetilde{u}_{c}|^2dx+\int_{\mathbb{R}^3}|(\nabla-iA)(w_{n}-\widetilde{u}_{c})|^2dx+o_{n}(1),\\
\int_{\mathbb{R}^3}V_{\Omega^{*}}(x)|w_{n}|^2dx&=\int_{\mathbb{R}^3}V_{\Omega^{*}}(x)|\widetilde{u}_{c}|^2dx
+\int_{\mathbb{R}^3}V_{\Omega^{*}}(x)|w_{n}-\widetilde{u}_{c}|^2dx+o_{n}(1),
\end{align*}
which directly follows \cite[(3.8) and (3.9)]{Dinh2022}. Hence,
\begin{align}\label{WE}
E_{\Omega^{*}}(w_{n})=E_{\Omega^{*}}(\widetilde{u}_{c})+E_{\Omega^{*}}(w_{n}-\widetilde{u}_{c})+o_{n}(1).
\end{align}
Next, we can show that for any $0<c<\kappa<\min\{c_0(\rho), c_1(\rho)\}$,
\begin{align}\label{SA}
\kappa m_{\Omega^{*}}(c,\rho)>cm_{\Omega^{*}}(\kappa,\rho),
\end{align}
where $m_{\Omega^{*}}(c,\rho)$ is defined in \eqref{4.4}. In fact, let $\{v_{n}\}\subseteq S(c)\cap B_{\Omega^{*}}(\rho)$ be such that
$\lim\limits_{n\to\infty}E_{\Omega^{*}}(v_{n})= m_{\Omega^{*}}(c,\rho)$. From \eqref{V1} and \eqref{V2},
we can assume $v_{n}\in B_{\Omega^{*}}(c\rho)$ for any $n$ large enough when $c<c_1(\rho)$. We further deduce that
\begin{align*}
\sqrt{\frac{\kappa}{c}}v_{n}\in S(\kappa)\cap B_{\Omega^{*}}(\kappa\rho)\subseteq S(\kappa)\cap B_{\Omega^{*}}(\rho),
\end{align*}
where $\kappa<1$. By Lemma \ref{lem16} and $N(v_{n})<0$, we derive that
\begin{align*}
m_{\Omega^{*}}(\kappa,\rho)&\leq E_{\Omega^{*}}(\sqrt{\frac{\kappa}{c}}v_{n})\\
&=\frac{\kappa}{2c}\int_{\mathbb{R}^3}|(\nabla -iA)v_{n}|^2+\frac{\kappa}{c}\int_{\mathbb{R}^3}V_{\Omega^{*}}|v_{n}|^2dx+\frac{\kappa^2}{2c^2}N(v_{n})\\
&=\frac{\kappa}{c}E_{\Omega^{*}}(v_{n})-\frac{\frac{\kappa}{c}(1-\frac{\kappa}{c})}{2}N(v_{n})\\
&=\frac{\kappa}{c}E_{\Omega^{*}}(v_{n})+\frac{\frac{\kappa}{c}(1-\frac{\kappa}{c})}{2}|N(v_{n})|\\
&<\frac{\kappa}{c}m_{\Omega^{*}}(c,\rho),
\end{align*}
which yields that \eqref{SA} holds. Moreover, setting $r_{n}:=|w_{n}-\widetilde{u}_{c}|_2^2$, we know that $r_{n}+r=c+o_{n}(1)$.
Then we can assume that $r_{n}\to l$ as $n\to \infty$, which states that $l+r=c$.
 Thus, in view of \eqref{WE} and \eqref{SA}, we conclude that
\begin{align*}
m_{\Omega^{*}}(c,\rho)=E_{\Omega^{*}}(w_{n})&=E_{\Omega^{*}}(\widetilde{u}_{c})+E_{\Omega^{*}}(w_{n}-\widetilde{u}_{c})+o_{n}(1)\\
&\geq m_{\Omega^{*}}(l, \rho)+m_{\Omega^{*}}(r, \rho)\\
&>\frac{l}{c}m_{\Omega^{*}}(c,\rho)+\frac{r}{c}m_{\Omega^{*}}(c,\rho)\\
&=m_{\Omega^{*}}(c,\rho),
\end{align*}
a contradiction. Thus, $|\widetilde{u}_{c}|_2^2=c$ and $w_{n}\to \widetilde{u}_{c}$ in $L^2(\mathbb{R}^3)$.
We further infer that $w_{n}\to \widetilde{u}_{c}$ in $L^4(\mathbb{R}^3)$ by the interpolation inequality. By Lemma \ref{lem2}, one has $N(w_{n})=N(\widetilde{u}_{c})+o_{n}(1)$.
It follows from the weak lower semi-continuity of norm that
\begin{align*}
m_{\Omega^{*}}(c,\rho)&=\liminf\limits_{n\to\infty}E_{\Omega^{*}}(w_{n})\\
&=\liminf\limits_{n\to\infty}\left[\int_{\mathbb{R}^3}|(\nabla-iA)w_{n}|^2dx+\int_{\mathbb{R}^3}V_{\Omega^{*}}(x)|w_{n}|^2dx+\frac{1}{2}N(w_{n})\right]\\
&\geq\int_{\mathbb{R}^3}|(\nabla-iA)\widetilde{u}_{c}|^2dx+\int_{\mathbb{R}^3}V_{\Omega^{*}}(x)|\widetilde{u}_{c}|^2dx+\frac{1}{2}N(\widetilde{u}_{c})\\
&=E_{\Omega^{*}}(\widetilde{u}_{c})\geq m_{\Omega^{*}}(c,\rho).
\end{align*}
It yields that $E_{\Omega^{*}}(\widetilde{u}_{c})= m_{\Omega^{*}}(c,\rho)$. Therefore, $\widetilde u_c$ is a minimizer for $m_{\Omega^{*}}(c,\rho)$ and $w_{n}\to \widetilde{u}_{c}$ in $\Sigma_{\Omega^{*}}$.

In what follows, we prove that $\widetilde{u}_{c}\not\in \partial B_{\Omega^{*}}(\rho)$. Obviously,
it suffices to verify that
$$
\mathcal{M}_{\Omega^{*}}^{\rho}(c)\subseteq B_{\Omega^{*}}(\rho c),
$$
where
$$\mathcal{M}_{\Omega^{*}}^{\rho}(c):=\{u\in S(c)\cap B_{\Omega^{*}}(\rho): E_{\Omega^{*}}(u)=m_{\Omega^{*}}(c,\rho)\}.
$$
Assume by contradiction that there exists $\psi\notin B_{\Omega^{*}}(\rho c)$ but $\psi\in\mathcal{M}_{\Omega^{*}}^{\rho}(c)$. It follows from \eqref{V2} that
\begin{align*}
m_{\Omega^{*}}(c,\rho)\leq&\inf\limits_{u\in S(c)\cap B_{\Omega^{*}}(\rho\frac{c}{2})}E_{\Omega^{*}}(u)\\
<&\inf\limits_{u\in S(c)\cap \left(B_{\Omega^{*}}(\rho)\setminus B_{\Omega^{*}}(\rho c)\right)}E_{\Omega^{*}}(u)\\
\leq& E_{\Omega^{*}}(\psi)=m_{\Omega^{*}}(c,\rho),
\end{align*}
which is a contradiction. Thus, $\widetilde{u}_{c}$ is a critical point of $E_{\Omega^{*}}|_{S(c)}$,
and there is a corresponding Lagrange multiplier $\widetilde{\mu}_{c}\in\mathbb{R}$
such that $E_{\Omega^{*}}'(\widetilde{u}_{c})+\widetilde{\mu}_{c}\widetilde{u}_{c}=0$ in $(\Sigma_{\Omega^{*}})^{-1}$.
Then, by Lemma \ref{lem2} and \eqref{DI} as well as \eqref{GN}, we deduce that
\begin{align*}
\widetilde{\mu}_{c}|\widetilde{u}_{c}|_2^2=&-\frac{1}{2}\int_{\mathbb{R}^3}|(\nabla-iA)\widetilde{u}_{c}|^2dx
-\int_{\mathbb{R}^3}V_{\Omega^{*}}(x)|\widetilde{u}_{c}|^2dx-N(\widetilde{u}_{c})\\
\leq&-\frac{1}{2}\int_{\mathbb{R}^3}|(\nabla-iA)\widetilde{u}_{c}|^2dx
-\int_{\mathbb{R}^3}V_{\Omega^{*}}(x)|\widetilde{u}_{c}|^2dx+\Lambda\int_{\mathbb{R}^3}|\widetilde{u}_{c}|^4dx\\
\leq&-\frac{1}{2}\int_{\mathbb{R}^3}|(\nabla-iA)\widetilde{u}_{c}|^2dx-\int_{\mathbb{R}^3}V_{\Omega^{*}}(x)|\widetilde{u}_{c}|^2dx
+\Lambda \mathcal C_4^4\left(\int_{\mathbb{R}^3}|\nabla|\widetilde{u}_{c}||^2dx\right)^{\frac{3}{2}}c^{\frac{1}{2}}\\
\leq&-\frac{1}{2}\int_{\mathbb{R}^3}|(\nabla-iA)\widetilde{u}_{c}|^2dx-\int_{\mathbb{R}^3}V_{\Omega^{*}}(x)|\widetilde{u}_{c}|^2dx
+\Lambda \mathcal C_4^4\left(\int_{\mathbb{R}^3}|(\nabla-iA)\widetilde{u}_{c}|^2dx\right)^{\frac{3}{2}}c^{\frac{1}{2}}\\
\leq&\|\widetilde{u}_{c}\|_{\Sigma_{\Omega^{*}}}^2\left(-\frac{1}{2}+\Lambda \mathcal C_4^4\rho^{\frac{1}{2}}c^{\frac{1}{2}}\right).
\end{align*}
Further, we choose $0<c<\bar c(\rho)$, where
\begin{align}\label{bar-c-rho}
\bar c(\rho):=\min\{c_0(\rho),c_1(\rho),\frac{1}{4\Lambda^2\mathcal C_4^8\rho}\}.
\end{align}
Then
$$
-\frac{1}{2}+\Lambda \mathcal C_4^4\rho^{\frac{1}{2}}c^{\frac{1}{2}}<0,
$$
which yields that $\widetilde{\mu}_{c}<0$. By \eqref{S1} and \eqref{S2}, we obtain that
\begin{align}\label{Z1}
\widetilde{\mu}_{c}\leq\zeta_{\Omega^{*}}^{0}\left(-\frac{1}{2}+\Lambda \mathcal C_4^4\rho^{\frac{1}{2}}c^{\frac{1}{2}}\right)
\leq\left(\gamma_2+\frac{\gamma_3}{2}\right)\Big(-\frac{1}{2}+\Lambda \mathcal C_4^4\rho^{\frac{1}{2}}c^{\frac{1}{2}}\Big).
\end{align}
On the other hand, since $(\lambda_1, \lambda_2)\in\mathbb{R}^2$ satisfies \eqref{UR'}, we have
\begin{align*}
\widetilde{\mu}_{c}|\widetilde{u}_{c}|_2^2=&-\frac{1}{2}\int_{\mathbb{R}^3}|(\nabla-iA)\widetilde{u}_{c}|^2dx
-\int_{\mathbb{R}^3}V_{\Omega^{*}}(x)|\widetilde{u}_{c}|^2dx-N(\widetilde{u}_{c})\\
=&-E_{\Omega^{*}}(\widetilde{u}_{c})-\frac{1}{2}N(\widetilde{u}_{c})>-E_{\Omega^{*}}(\widetilde{u}_{c})=-m_{\Omega^{*}}(c,\rho).
\end{align*}
From Lemma \ref{lem16}, one has
\begin{align}\label{Z2}
\widetilde{\mu}_{c}>-\frac{m_{\Omega^{*}}(c,\rho)}{c}\geq-\left(\gamma_2+\frac{\gamma_3}{2}\right).
\end{align}
Combining \eqref{Z1} and \eqref{Z2}, we obtain \eqref{T3}.

$(ii)$ Let $\widetilde{u}_{c}\in \mathcal{M}_{\Omega^{*}}^{\rho}(c)$ is a minimizer for $m_{\Omega^{*}}(c,\rho)$ obtained by $(i)$.
From Lemma \ref{lem2}, \eqref{GN} and \eqref{DI}, we infer that
\begin{align*}
\int_{\mathbb{R}^3}&|(\nabla-iA)\widetilde{u}_{c}|^2+2V_{\Omega^{*}}(x)|\widetilde{u}_{c}|^2dx\\
&=2E_{\Omega^{*}}(\widetilde{u}_{c})-N(\widetilde{u}_{c})\\
&\leq2m_{\Omega^{*}}(c,\rho)+\Lambda\int_{\mathbb{R}^3}|\widetilde{u}_{c}|^4dx\\
&\leq2m_{\Omega^{*}}(c,\rho)+\Lambda \mathcal C_4^4\left(\int_{\mathbb{R}^3}|\nabla|\widetilde{u}_{c}||^2dx\right)^{\frac{3}{2}}c^{\frac{1}{2}}\\
&\leq2m_{\Omega^{*}}(c,\rho)++\Lambda \mathcal C_4^4\left(\int_{\mathbb{R}^3}|(\nabla-iA)\widetilde{u}_{c}|^2dx\right)^{\frac{3}{2}}c^{\frac{1}{2}}\\
&<2(\gamma_2+\frac{\gamma_3}{2})c+\Lambda \mathcal C_4^4\rho^{\frac{3}{2}}c^{\frac{1}{2}},
\end{align*}
which means that
\begin{align*}
\|\widetilde{u}_{c}\|_{\Sigma_{\Omega^{*}}}^2<2(\gamma_2+\frac{\gamma_3}{2}+\frac{1}{2})c+\Lambda \mathcal C_4^4\rho^{\frac{1}{2}}c^{\frac{1}{2}}.
\end{align*}
Thus $\widetilde{u}_{c}\to 0$ in $\Sigma_{\Omega^{*}}$ as $c\to0^{+}$.
Next we recall from \eqref{T3} that
$$-(\gamma_2+\frac{\gamma_3}{2})<\widetilde{\mu}_{c}\leq(\gamma_2+\frac{\gamma_3}{2})\Big(-\frac{1}{2}+\Lambda \mathcal C_4^4\rho^{\frac{1}{2}}c^{\frac{1}{2}}\Big)
$$
for all $0<c<\bar c(\rho)$. Then there
exists a constant $\widetilde{\mu}_0\in[-(\gamma_2+\frac{\gamma_3}{2}),-\frac{1}{2}(\gamma_2+\frac{\gamma_3}{2})] $
such that $\widetilde{\mu}_{c}\to\widetilde{\mu}_0 $ as $c\to0^{+}$.
By Lemma \ref{lem2}, \eqref{GN}, \eqref{DI} and the fact $\widetilde{u}_{c}\in B_{\Omega^{*}}(c\rho)$, we can get that
\begin{align*}
0<\frac{|N(\widetilde{u}_{c})|}{|\widetilde{u}_{c}|_2^2}\leq\frac{\Lambda |\widetilde{u}_{c}|_4^4}{|\widetilde{u}_{c}|_2^2}\leq\frac{\Lambda \mathcal C_4^4|\nabla |\widetilde{u}_{c}|_2^3}{|\widetilde{u}_{c}|_2}\leq\frac{\Lambda \mathcal C_4^4|(\nabla-iA)\widetilde{u}_{c} |_2^3}{|\widetilde{u}_{c}|_2}\leq\Lambda \mathcal C_4^4\rho^{\frac{3}{2}}c\to0
\end{align*}
as $c\to0^{+}$.
It results that
\begin{align*}
\lim\limits_{c\to0^{+}}\frac{m_{\Omega^{*}}(c,\rho)}{c}
=&\lim\limits_{c\to0^{+}}\left[\frac{1}{|\widetilde{u}_{c}|_2^2}\int_{\mathbb{R}^3}|(\nabla-iA)\widetilde{u}_{c}|^2dx
+\frac{1}{|\widetilde{u}_{c}|_2^2}\int_{\mathbb{R}^3}V_{\Omega^{*}}(x)\widetilde{u}_{c}|^2dx
+\frac{1}{2}\frac{N(\widetilde{u}_{c})}{|\widetilde{u}_{c}|_2^2}\right]\\
=&\lim\limits_{c\to0^{+}}\frac{1}{|\widetilde{u}_{c}|_2^2}\left[\frac{1}{2}\int_{\mathbb{R}^3}|(\nabla-iA)\widetilde{u}_{c}|^2dx
+\int_{\mathbb{R}^3}V_{\Omega^{*}}(x)|\widetilde{u}_{c}|^2dx\right].
\end{align*}
Similarly, we have
\begin{align*}
\lim\limits_{c\to0^{+}}\widetilde{\mu}_{c}=&\lim\limits_{c\to0^{+}}\left[-\frac{1}{2|\widetilde{u}_{c}|_2^2}\int_{\mathbb{R}^3}|(\nabla-iA)\widetilde{u}_{c}|^2dx
-\frac{1}{|\widetilde{u}_{c}|_2^2}\int_{\mathbb{R}^3}V_{\Omega^{*}}(x)|\widetilde{u}_{c}|^2dx-\frac{N(\widetilde{u}_{c})}{|\widetilde{u}_{c}|_2^2}\right]\\
=&\lim\limits_{c\to0^{+}}\frac{1}{|\widetilde{u}_{c}|_2^2}\left[-\frac{1}{2}\int_{\mathbb{R}^3}|(\nabla-iA)\widetilde{u}_{c}|^2dx
-\int_{\mathbb{R}^3}V_{\Omega^{*}}(x)|\widetilde{u}_{c}|^2dx\right].
\end{align*}
Thus,
$$
\lim\limits_{c\to0^{+}}\frac{m_{\Omega^{*}}(c,\rho)}{c}
=\lim\limits_{c\to0^{+}}-\widetilde{\mu}_{c}=-\widetilde{\mu}_0.
$$ This completes the proof.
\qed
\\~
\ \\
\noindent{\bf Conflict of interest:} The authors declare that they have no conflict of interest.\\
\ \\
\noindent{\bf Acknowledgments:} This research is supported by National Natural Science Foundation of China (No. 12371120).


\begin{thebibliography}{99}

\bibitem{Antoine2018}
\newblock X. Antoine, Q. Tang, Y. Zhang,
\newblock A preconditioned conjugated gradient method for computing ground states of rotating
dipolar Bose-Einstein condensates via kernel truncation method for dipole-dipole interaction evaluation,
\newblock \emph{Commun. Comput. Phys.}, \textbf{24} (2018), 966--988.

\bibitem{Antonelli2015}
\newblock P. Antonelli, R. Carles, J.D. Silva,
\newblock Scattering for nonlinear Schr\"{o}dinger equation under partial harmonic confinement,
\newblock \emph{Commun. Math. Phys.}, \textbf{334} (2015), 367--396.

\bibitem{Antonelli2011}
\newblock P. Antonelli, C. Sparber,
\newblock Existence of solitary waves in dipolar quantum gases,
\newblock \emph{Phys. D}, \textbf{240} (2011), 426--431.

\bibitem{Arbunich2019}
\newblock J. Arbunich, I. Nenciu, C. Sparber,
\newblock Stability and instability properties of rotating Bose-Einstein condensates,
\newblock \emph{Lett. Math. Phys.}, \textbf{109} (2019), 1415--1432.

\bibitem{Jeanjean2017}
\newblock J. Bellazzini, L. Jeanjean,
\newblock On dipolar quantum gases in the unstable regime,
\newblock \emph{SIAM J. Math. Anal.}, \textbf{48} (2016), 2028--2058.

\bibitem{Bellazzini2017}
\newblock J. Bellazzini, N. Boussa\"{i}d, L. Jeanjean, et al.,
\newblock Existence and stability of standing waves for supercritical NLS with a partial confinement,
\newblock \emph{Commun. Math. Phys.}, \textbf{353} (2017), 229--251.

\bibitem{Bao2005}
\newblock W. Bao, H. Wang, P.A. Markowich,
\newblock Ground, symmetric and central vortex states in rotating Bose-Einstein condensates,
\newblock \emph{Commun. Math. Sci.}, \textbf{3} (2005), 57--88.

\bibitem{Bao2010}
\newblock W. Bao, Y. Cai, H. Wang,
\newblock Efficient numerical methods for computing ground states and dynamics of dipolar Bose-Einstein condensates,
\newblock \emph{J. Comput. Phys.}, \textbf{229} (2010), 7874-7892.


\bibitem{Bao2012}
\newblock W. Bao, B.A. Naoufel, Y. Cai,
\newblock Gross-Pitaevskii-Poisson equations for dipolar Bose-Einstein condensate with anisotropic confinement,
\newblock \emph{SIAM J. Math. Anal.}, \textbf{44} (2012), 1713-1741.

\bibitem{CaiYY2011}
\newblock Y.Y. Cai,
\newblock Mathematical theory and numerical methods for the Gross-Piatevskii equations and applications,
\newblock PH.D. Thesis, National University of Singapore (2011).



\bibitem{Carles2015}
\newblock R. Carles, H. Hajaiej,
\newblock Complementary study of the standing wave solutions of the Gross-Pitaevskii equation in dipolar quantum gases,
\newblock \emph{Bull. London Math. Soc.}, \textbf{47} (2015), 509--518.

\bibitem{Carles2008}
\newblock R. Carles, P. Markowich, C. Sparber,
\newblock On the Gross-Pitaevskii equation for trpped dipolar quantum gases,
\newblock \emph{Nonlinearity}, \textbf{21} (2008), 2569--2590.

\bibitem{Dell2007}
\newblock D.H.J. O'Dell, C. Eberlein,
\newblock Vortex in a trapped Bose-Einstein condensate with dipole-dipole interactions,
\newblock \emph{Phys. Rev. A}, \textbf{75} (2007), 013604.

\bibitem{Dinh2021}
\newblock V.D. Dinh,
\newblock On the instability of standing waves for 3D dipolar Bose-Einstein condensates,
\newblock \emph{Phys. D}, \textbf{419} (2021), 12pp.

\bibitem{Dinh2022}
\newblock V.D. Dinh,
\newblock Existence and stability of standing waves for nonlinear Schr\"{o}dinger equations with a critical rotational speed,
\newblock \emph{Lett. Math. Phys.}, \textbf{112} (2022), 36pp.

\bibitem{Esteban1989}
\newblock M.J. Esteban, P.L. Lions,
\newblock Stationary solutions of nonlinear Schr\"{o}dinger equations with an external magnetic field. In: Partial Differential Equations and the Calculus of Variations,
\newblock Vol. I, Progr. Nonlinear Differential Equations Appl., 1, Birkh\"{a}user Boston, Boston, MA, 1989.

\bibitem{Feng2021}
\newblock B. Feng, Q. Wang,
\newblock Strong instability of standing waves for the nonlinear Schr\"{o}dinger equation in trapped dipolar quantum gases,
\newblock \emph{J. Dynam. Differential Equations}, \textbf{33} (2021) 1989--2008.


\bibitem{Jeanjean1997}
\newblock L. Jeanjean,
\newblock Existence of solutions with prescribed norm for semilinear elliptic equations,
\newblock \emph{Nonlinear Anal.}, \textbf{28} (1997), 1633--1659.

\bibitem{Lewin2018}
\newblock M. Lewin, P.T. Nam, N. Rougerie,
\newblock Blow-up profile of rotating 2D focusing Bose gases,
\newblock Springer Verlag, Macroscopic Limits of Quantum Systems 2018.

\bibitem{Lieb2001}
\newblock E.H. Lieb, M. Loss,
\newblock Analysis, Graduate Studies in Mathematics,
\newblock Vol. 14, 2nd edn. American Mathematical Society, Providence (2001).

\bibitem{Liu2015}
\newblock B. Liu, L. Ma,
\newblock Blow up threshold for the Gross-Pitaevskii system with combined nonlocal nonlinearities,
\newblock \emph{J. Math. Anal. Appl.}, \textbf{425} (2015), 1214--1224.

\bibitem{Luo2021}
\newblock X. Luo, T. Yang,
\newblock Multiplicity, asymptotics and stability of standing waves for nonlinear Schr\"{o}dinger equation with rotation,
\newblock \emph{J. Differential Equations}, \textbf{304} (2021), 326--347.

\bibitem{Ma2011}
\newblock L. Ma, P. Cao,
\newblock The threshold for the focusing Gross-Pitaevskii equation with trapped dipolar quantum gases,
\newblock \emph{J. Math. Anal. Appl.}, \textbf{381} (2011), 240--246.

\bibitem{Ma2013}
\newblock L. Ma, J. Wang,
\newblock Sharp threshold of the Gross-Pitaevskii equation with trapped dipolar quantum gases,
\newblock \emph{Canad. Math. Bull.}, \textbf{56} (2013), 378--387.

\bibitem{Tang2017}
\newblock Q. Tang, Y. Zhang, N.J. Mauser,
\newblock A robust and efficient numerical method to compute the dynamics of the rotating two-component dipolar Bose-Einstein condensates,
\newblock \emph{Comput. Phys. Commun.}, \textbf{219} (2017), 223--235.

\bibitem{Triay2018}
\newblock A. Triay,
\newblock Derivation of the dipolar Gross-Pitaevskii energy,
\newblock \emph{SIAM J. Math. Anal.}, \textbf{50} (2018), 33--63.

\bibitem{Weinstein2009}
\newblock M. Weinstein,
\newblock Nonlinear Schr\"{o}dinger equations and sharp interpolation estimates,
\newblock \emph{Commun. Math. Phy.}, \textbf{87} (1983), 567--576.

\bibitem{Minimax}
\newblock M. Willem,
\newblock Minimax theorems, vol.~24 of Progress in Nonlinear
  Differential Equations and their Applications,
\newblock Birkh\"{a}user Boston, Inc., Boston, MA, 1996.

\bibitem{Wu2024}
\newblock M.-H. Wu, C.-L. Tang,
\newblock Multiplicity of normalized solutions for dipolar Gross-Pitaevskii equation with a mass subcritical perturbation,
\newblock \emph{J. Fixed Point Theory Appl}, \textbf{26} (2024), 25pp.

\bibitem{Yi2006}
\newblock S. Yi, H. Pu,
\newblock Vortex structres in dipolar condensates,
\newblock \emph{Phys. Rev. A}, \textbf{73} (2006), 061602.

\bibitem{Zhang2000}
\newblock J. Zhang,
\newblock Stability of standing waves for nonlinear Schr\"{o}dinger equations with unbounded potentials,
\newblock \emph{Z. Angew. Math. Phys.}, \textbf{51} (2000), 498--503.

\bibitem{Zhang2016}
\newblock X.-F. Zhang, W. Han, H.-F. Jiang, et al.,
\newblock Topological defect formation in rotating binary dipolar Bose-Einstein condensate,
\newblock \emph{Ann. Physics}, \textbf{375} (2016), 368--377.

\end{thebibliography}


\end{document}